\documentclass[]{elsarticle}

\journal{Journal of Computational Physics}

\usepackage{hyperref}

\usepackage{booktabs}
\usepackage{mathtools}
\usepackage[ruled,vlined]{algorithm2e}
\usepackage{subfig}

\usepackage{amsmath}
\usepackage{amssymb,amsfonts}
\usepackage{amsthm}
\newtheorem{remark}{Remark}
\newtheorem{definition}{Definition}
\newtheorem{proposition}{Proposition}

\usepackage{microtype}

\usepackage[american]{babel}
\usepackage[utf8]{inputenc}
\usepackage[T1]{fontenc}

\usepackage{lmodern}

\newcommand{\R}{\mathbb{R}}
\newcommand{\C}{\mathbb{C}}
\newcommand{\e}{\varepsilon}

\newcommand{\er}{\tilde\varepsilon_{r}\,}
\newcommand{\mur}{\mu_{r}\,}
\newcommand{\muri}{\mu_{r}^{-1}\,}
\newcommand{\ssr}{\sigma^\Sigma_r}
\newcommand{\tissr}{{\breve\sigma}^\Sigma_r}
\newcommand{\km}{k_{m}}
\newcommand{\kmr}{k_{m,r}}

\renewcommand{\vec}[1]{\boldsymbol{#1}}

\newcommand{\avE}{\vec {\acute E}}
\newcommand{\hvB}{\vec {\hat B}}
\newcommand{\hvE}{\vec {\hat E}}
\newcommand{\hvJa}{\vec {\hat J}_{a}}
\newcommand{\hvMa}{\vec {\hat M}_a}
\newcommand{\J}{\mathcal{J}}
\newcommand{\vB}{\vec B}
\newcommand{\vE}{\vec E}
\newcommand{\ve}{\vec e}
\newcommand{\vF}{\vec F}
\newcommand{\vJa}{\vec J_{a}}
\newcommand{\vMa}{\vec M_a}
\newcommand{\vMF}{\vec{\mathcal{F}}}
\newcommand{\vME}{\vec{\mathcal{E}}}
\newcommand{\vMB}{\vec{\mathcal{B}}}
\newcommand{\vMMa}{\vec{\mathcal{M}}_a}
\newcommand{\vMJa}{\vec{\mathcal{J}}_a}
\newcommand{\vn}{\vec \nu}
\newcommand{\vp}{\vec \varphi}
\newcommand{\vx}{\vec x}
\newcommand{\vZ}{\vec Z}
\newcommand{\vU}{\vec U}

\newcommand{\TTH}{\mathbb{T}_H}

\newcommand{\dx}{\,{\mathrm d}x}
\newcommand{\dox}{\,{\mathrm d}o_x}

\begin{document}


\begin{frontmatter}

  \title{%
    Dipole excitation of surface plasmon on a conducting sheet:
    finite element approximation and validation
  }

  \author[minneapolis,murisupport]{Matthias Maier}
  \ead{msmaier@umn.edu}

  \author[maryland1,nsfsupport]{Dionisios Margetis}
  \ead{dio@math.umd.edu}

  \author[minneapolis,murisupport]{Mitchell Luskin\corref{corresponding}}
  \cortext[corresponding]{Corresponding author}
  \ead{luskin@umn.edu}

  \fntext[murisupport]{%
    The first and third authors (MM and ML) were supported in part by ARO
    MURI Award W911NF-14-1-0247. The third author (ML) was also supported
    in part by the Radcliffe Institute for Advanced Study at Harvard
    University.}

  \fntext[nsfsupport]{This author’s research was supported in part by NSF
    DMS-1412769.}

  \address[minneapolis]{School of Mathematics, University of Minnesota,
    Minneapolis, Minnesota 55455, USA.}
  \address[maryland1]{Department of Mathematics, and Institute for Physical
    Science and Technology, and Center for Scientific Computation and
    Mathematical Modeling, University of Maryland, College Park, Maryland
    20742, USA.}

  \begin{abstract}
    We formulate and validate a finite element approach to the propagation of
    a slowly decaying electromagnetic wave, called surface plasmon-polariton,
    excited along a conducting sheet, e.g., a single-layer graphene sheet, by
    an electric Hertzian dipole. By using a suitably rescaled form of
    time-harmonic Maxwell's equations, we derive a variational formulation
    that enables a direct numerical treatment of the associated class of
    boundary value problems by appropriate curl-conforming finite elements.
    The conducting sheet is modeled as an idealized hypersurface with an
    effective electric conductivity. The requisite weak discontinuity for the
    tangential magnetic field across the hypersurface can be incorporated
    naturally into the variational formulation. We carry out numerical
    simulations for an infinite sheet with constant isotropic conductivity
    embedded in two spatial dimensions; and validate our numerics against the
    closed-form exact solution obtained by the Fourier transform in the
    tangential coordinate. Numerical aspects of our treatment such as an
    absorbing perfectly matched layer, as well as local refinement and
    a posteriori error control are discussed.
  \end{abstract}


\end{frontmatter}



\section{Introduction}
\label{sec:Intro}

The manipulation of the electronic structure of low-dimensional materials
has recently been the subject of active research with applications in
spintronics, quantum information processing, and energy
storage~\cite{castroneto09,pesin12,trauzettel07,green12,geim13,torres-book,zhang-book}.
In particular, the electric conductivity of atomically thick materials such
as graphene and black phosphorous yields an effective complex permittivity
with a {\em negative}  real part in the infrared
spectrum~\cite{zhang-book,low14}. This feature allows for the propagation
of slowly decaying electromagnetic waves, called {\em surface
plasmons-polaritons} (SPPs), that are confined near the material interface
with wavelengths much shorter than the wavelength of the free-space
radiation~\cite{pitarke07,zhang-book}. These SPPs are promising ingredients
in the design of ultrafast photonic circuits~\cite{pitarke07}.

Experimental efforts to generate SPPs focus on the requisite phase matching
between waves sustained in free space and the material of
interest~\cite{bludov13}. This matching is enabled by the excitation of
sufficiently large wave numbers tangential to the associated interface. A
technique to achieve this goal is to place a resonant antenna on a graphene
sheet~\cite{gonzalez14,liu12}. The prediction of the resulting waves relies
on solving an intricate boundary value problem for the time-harmonic
Maxwell equations.

Our purpose with this paper is three-fold. First, we aim to develop a
general variational  framework for the numerical treatment of
electromagnetic-wave propagation in conducting materials modeled as
hypersurfaces in a Euclidean space of arbitrary dimension. These
hypersurfaces are irradiated by fields produced by external, compactly
supported current-carrying sources operating at a fixed, yet arbitrary,
frequency. Second, we  validate this framework by comparison of finite
element numerical simulations to an exact solution of Maxwell's equations
in the simplified case when the conducting sheet is infinite and embedded
in a two-dimensional (2D) space in the presence of a Hertzian electric
dipole directed vertically to the sheet. Third, we demonstrate that our
numerical approach is capable of efficiently resolving SPP structures. To
this end, goal-oriented adaptivity for local mesh refinement and a
perfectly matched layer are tailored to the fine structure of SPPs.

This work has been broadly motivated by the growing urge to engineer
microscopic details of low-dimensional conducting materials in order to
establish desired optical properties at larger scales~\cite{pitarke07}. A
key objective is to elucidate how the electric conductivity of the material
affects characteristics of electromagnetic wave propagation~\cite{cheng13}.
There is a compelling need for controllable numerical schemes which, placed
on firm mathematical grounds, can reliably describe the SPPs in a variety
of geometries. Our paper offers a systematic approach to solving this
problem. In particular:
\begin{itemize}
\item
  We formulate a variational framework suitable for the finite element
  treatment of electromagnetic wave propagation along conducting sheets
  embedded in spaces of arbitrary dimensions in the presence of external
  sources (Section~\ref{subsec:variational}).
\item
  We discretize our variational formulation by use of appropriate
  curl-conforming Nédélec-elements. This scheme is implemented in a modern
  C++ framework~\cite{dealii83} and accounts for the fine scale of the SPP
  as well as for the requisite radiation condition at infinity via a {\em
  perfectly matched layer} (PML) (Section~\ref{sec:numerics}). Our approach
  does not require regularization of the conducting sheet.
\item
  We validate our numerical treatment by comparison of numerics to a
  tractable exact solution in the simplified case with a vertical electric
  dipole radiating over an infinite, planar conducting sheet; for
  convenience, we consider an ambient 2D Euclidean space. In particular, we
  numerically single out the SPP; and also compare it to the slowly-varying
  radiation field (see Sections~\ref{sec:analytical}
  and~\ref{sec:computationalresults}).
\end{itemize}

The extensive literature in plasmonics attests to the rich variety of
computational methods and tools; see, e.g.,~\cite{gallinet2015,yeh-book}.
For example, an approach is to model graphene as a region with finite
thickness (as opposed to a boundary)~\cite{koppens2011}. In the present
work, we demonstrate the ability of curl-conforming Nédélec-elements to
accurately capture the fine scale of the SPPs, by replacing the conducting
sheet with a set of boundary conditions on a hypersurface.

Note that our setting and objectives in this work, focusing on the
propagation of SPPs along low-dimensional conducting materials, are
distinctly different from the modeling and computation of the interaction
of plasmonic nanoparticles with electromagnetic
fields~\cite{chau2011,edel2016}.  Our focus on the computation of SPPs
along low-dimensional conducting materials and its validation against
recently derived analytic solutions~\cite{margetis15} is also distinct from
the more classical study of surface plasmons on bulk
materials~\cite{raether86}.  However, our formulation and general approach
is also applicable in this setting.

Throughout the paper, we assume that the reader is familiar with the
fundamentals of classical electromagnetic wave theory; for extensive and
comprehensive treatments of this subject, see, e.g.,~\cite{king92,
muller69, schwartz72}.

\subsection{Motivation: Surface plasmonics}
\label{subsec:SPP-intro}

Atomically thick conducting materials such as graphene, black phosphorus,
and van der Waals heterostructures have been the focus of intensive
studies~\cite{geim13,pitarke07}. The dispersion
relations of these structures for electromagnetic wave propagation have novel features. The
implications of this dispersion at the infrared spectrum is a
theme of essence in surface plasmonics~\cite{samaier07,zhang-book}. Specifically, in
the terahertz frequency range, the effective electric conductivity,
$\sigma^\Sigma$, emerging from the coordinated motion of quasi-free
electrons, can have an appreciable imaginary part. Furthermore, it has been predicted via
numerical simulations that the decoration of graphene by chains of organic
molecules may result in a dramatic alteration of
$\sigma^\Sigma$~\cite{cheng13,dewapriya15}. This prediction paves the way to
unconventional means of controlling electronic transport.

From the viewpoint of Maxwell's equations, the effective dielectric
permittivity of a conducting sheet may have a negative real part. The
resulting {\em metamaterial} has optical properties different from those of
a conventional conductor~\cite{zhang-book}. In particular,
electromagnetic waves of transverse-magnetic (TM) polarization possibly
propagating through the atomically thick material are characterized by a
dispersion relation that allows for transmitted wave numbers much larger
than the free-space wave number, $k$. For an isotropic and homogeneous
ambient space, with wave number $k=\omega\sqrt{\tilde\e\mu}$ and scalar and
$\vx$-independent $\mu$ and $\tilde\e$, and an isotropic and homogeneous
conducting sheet, the condition $|\omega\mu\sigma^\Sigma|\ll |k|$ yields
the simplified dispersion relation~\cite{samaier07,hanson08,raether86}
\begin{equation}\label{eq:simple_disp-reln}
\sqrt{k^2-k_x^2}\approx -\left(\frac{2k}{\omega\mu\sigma^\Sigma}\right)\,k
\end{equation}
for TM waves; $k_x$ denotes the wave number tangential to the material
interface. Hence, if $k$ is positive, \eqref{eq:simple_disp-reln} has an
admissible solution, $k_x$, provided $\text{Im}\sigma^\Sigma >0$ under an
assumed $e^{-i\omega t}$ time dependence. Note that $|k_x|\gg |k|$.

The excitation of SPPs on a
homogeneous conducting sheet cannot be achieved by direct illumination of
the material by an incident plane wave.  There is intrinsic need for phase matching between the
waves propagating in different materials, e.g., air and conducting
sheet~\cite{bludov13,samaier07}. An indirect means of establishing this
matching is to add metal contacts to the
interface~\cite{satou07}.
More generally, it is plausible to prescribe current-carrying sources of
compact support that optimize attributes of the SPP by variation of the
frequency or size of source or the esheet conductivity~\cite{cheng13,liu12,gonzalez14}.

\subsection{Our approach}
\label{subsec:main_results}

A jump condition created by an electric conductivity on an interface is a
key ingredient in the modeling of thin conducting materials such as
graphene, black phosphorus, and a variety of heterostructures resulting
from stacking a few distinct crystalline sheets on top of each other. One
of our tasks with this work is to construct a variational formulation,
well-suited for the finite element method, that naturally incorporates such
a jump condition in the presence of external sources.

The variational formulation is implemented by utilizing an appropriate
curl-conforming finite element space that only enforces the continuity of
the tangential components across elements. The use of \emph{higher-order
conforming} elements is well suited for the numerical problem at hand. The
weak discontinuity across the interface can be aligned with the
triangulation and the regularity of the solution away from the interface
leads to high convergence rates. For overcoming the two-scale character
with much finer SPP structures close to the interface, an adaptive, local
refinement strategy based on a posteriori error estimates is used. The
a posteriori error estimates are computed by solving an adjoint problem
(\emph{dual weighted residual method}) \cite{becker1996a} and lead to
optimally refined meshes.

Notably, our approach does not require the regularization of the conducting
sheet by a layer with artificial thickness. Instead, the sheet can be
directly approximated as a lower-dimensional interface. Further, we treat
the full scattering problem with an incident wave generated by a Hertzian
dipole source instead of merely solving the associated eigenvalue problem
for the SPP.

For validation of our treatment, the finite element computations stemming
from our approach are compared to the exact solution of Maxwell's equations
for a vertical electric dipole over an isotropic and homogeneous conducting
sheet in 2D. In this case, all field components are expressed via 1D
Fourier integrals and, thus, are amenable to accurate numerical
integration. By this formalism, the SPP is defined as the contribution from
a simple pole in the Fourier domain. This contribution is to be contrasted
to the slowly-varying radiation field. Our numerics indicate that the SPP
dominates the scattered field at distances of the order of the free-space
wavelength from the dipole source, in agreement with analytical estimates from the exact solution.

\subsection{Related work}
\label{subsec:past}
Electromagnetic wave propagation along boundaries, especially
the boundary separating air and earth or sea, has been the subject of studies for over a
century. A review can be found in~\cite{king92}.

This insight is valuable yet insufficient for plasmonic applications
related to low-dimensional materials. It is compelling to consider
implications of the {\em metamaterial} character of {\em atomically} thick
conducting sheets in the terahertz frequency range~\cite{bludov13}. In
particular, in the presence of an electric Hertzian dipole source, boundary
condition~\eqref{eq:jumpcondition} with $\text{Im}\sigma^\Sigma >0$ can
result in a SPP~\cite{bludov13}, to be contrasted to surface waves in
radio-frequencies which have wave numbers nearly equal to the free-space
one~\cite{king92}.

In the last few decades, several groups have been studying
implications of surface plasmonics; for a (definitely non-exhaustive)
sample of related works,
see~\cite{bludov13,hanson08,hanson11,liu12,margetis15,nikitin11,raether86,samaier07,satou07}.
For instance, in~\cite{bludov13} the authors review macroscopic
properties of the electric conductivity of graphene, derive dispersion
relations for electromagnetic plane waves in inhomogeneous structures, and
discuss methods for exciting SPPs; see also the integral-equation
approach in~\cite{satou07}. On the other hand, the problem of a
radiating dipole source near a graphene sheet is semi-analytically
addressed in~\cite{hanson08,hanson11}.
In the same vein, in~\cite{nikitin11} the authors numerically
study the field produced by dipoles near a graphene sheet, recognizing a
region where the scattered field may be significant. Most recently, two of
us derived closed-form analytical expressions for the
electromagnetic field when the dipole source and observation point lie on
the sheet~\cite{margetis15}.

In the aforementioned works, the rigorous numerical treatment of Maxwell's
equations is {\em not} of primary concern. Our work here aims to build a
framework that places the finite element treatment of a variational
boundary value problem on firmer mathematical grounds. This opens up the
possibility of numerically studying SPPs in experimentally accessible
geometries in an error-controllable fashion.
The finite element treatment of Maxwell's equations is a well-established
area of research~\cite{brenner13, brenner16, monk03, nedelec86, nedelec01}.
In particular, our work is related to an adaptive finite element framework
for Maxwell's equations \cite{brenner13, brenner16} that uses a
residual-based a posteriori error estimator for adaptive local refinement
of a triangular mesh. Further, we point out that the well-posedness of our
variational approach (see Section~\ref{sec:variational}) is closely
connected to the question of well-posedness of time-harmonic Maxwell's
equations with sign-changing dielectric permittivity \cite{bendhia14}. The
reason for this connection is that the jump condition for a field component
across a boundary with a complex-valued conductivity (that we use to model
graphene) can be understood as the zero-thickness limit of a \emph{bulk
graphene} region with a negative dielectric permittivity.

\subsection{Open problems}
\label{subsec:limitations}
Our work here focuses on the development of a reasonably general
variational framework. We validate this formulation by comparison of the
ensuing finite element numerical computations to the exact solution of
Maxwell's equations in a relatively simple yet nontrivial geometry in 2D.
We deem these tasks as necessary first steps in establishing the proposed
numerical framework; these steps should precede applications to more
complicated cases of physical interest.

Therefore, our work admits several extensions and leaves a few pending
issues, from the viewpoints of both analysis and applications. For example,
we have not made attempts to fully characterize error estimates following
from our treatment. The subtleties related to the possible spatial
variation or anisotropy of surface conductivity, $\sigma^\Sigma$, lie
beyond our present scope. The numerics for a current-carrying source over a
conducting film in 3D~\cite{hanson08,margetis15} have not been carried out,
because of the expensive computations involved. The more elaborate yet
experimentally accessible case with a receiving antenna lying on the
material interface~\cite{gonzalez14}, where the current distribution on the
source forms part of the solution, is a promising topic of near-future
investigation.

\subsection{Outline of paper}
\label{subsec:outline}
The remainder of our paper is organized as follows. In
Section~\ref{sec:variational}, we provide the desired variational
characterization for boundary value
problem~\eqref{eq:timeharmonicmaxwell}--\eqref{eq:S-M_cond}.
Section~\ref{sec:analytical} focuses on the derivation of an exact solution
for~\eqref{eq:timeharmonicmaxwell}--\eqref{eq:S-M_cond} in 2D, assuming
that the external current-carrying source is a vertical electric dipole and
the conducting sheet is homogeneous and isotropic. In
Section~\ref{sec:numerics}, we describe the discretization of our
variational formulation in the context of finite elements; in particular,
we discuss the error control by our treatment
(Section~\ref{sec:aposteriori}). Section~\ref{sec:computationalresults}
present computational results stemming from our approach for an infinite
conducting sheet in 2D, along with comparisons with an exact solution.
Finally, Section~\ref{sec:conclusion} concludes our paper with a summary of
our results and an outlook.


\section{Variational formulation}
\label{sec:variational}

In this section, we derive a variational formulation for the time-harmonic
Maxwell equations with an interface jump condition. We introduce a slightly
modified rescaling of the associated equations to dimensionless forms that
are best suited for the numerical observation of the SPP in our treatment.
The interface jump condition (\ref{eq:jumpcondition}) enters the
variational formulation in the form of a \emph{weak discontinuity} (with
the second jump-condition for $\vE$ being naturally encoded in the {\em
ansatz} space).

\subsection{Preliminaries: Boundary value problem}
\label{subsec:Prelim}

Next, we formulate the corresponding boundary value problem for the
conducting sheet, emphasizing the discontinuity of the magnetic field
across the sheet. The starting point of our analysis is the strong form of
Maxwell's equations for the time-harmonic electromagnetic field,
$(\vME(\vx,t), \vMB(\vx,t))=\text{Re}\,\big\{e^{-i\omega t}(\vE(\vx),
\vB(\vx))\big\}$, viz.,~\cite{schwartz72}
\begin{align}
  \begin{cases}
    \begin{aligned}
      -i\omega\vB+\nabla\times\vE \;&=\; -\vMa~,
      \\[0.1em]
      \nabla\cdot\vB \;&=\; \frac 1{i\omega}
      \nabla\cdot\vMa~,
      \\[0.1em]
      i\omega\tilde\e\vE+\nabla\times\big(\mu^{-1}\vB\big) \;&=\; \vJa~,
      \\[0.1em]
      \nabla\cdot\big(\tilde\e\vE\big) \;&=\; \frac 1{i\omega}
      \nabla\cdot\vJa~.
    \end{aligned}
  \end{cases}
  \label{eq:timeharmonicmaxwell}
\end{align}
A few comments on~\eqref{eq:timeharmonicmaxwell} are in order. The
(constant) parameter $\omega$ is the temporal angular frequency ($\omega
>0$). We assume that all material parameters are time independent;
furthermore, the time-independent, externally applied electric- and
magnetic-current densities, $\vJa(\vx)$ and $\vMa(\vx)$, respectively,
arise from the time-harmonic densities
$\vMJa(\vx,t)=\text{Re}\,\big\{e^{-i\omega t}\vJa(\vx)\big\}$ and
$\vMMa(\vx,t)=\text{Re}\,\big\{e^{-i\omega t}\vMa(\vx)\big\}$. The
second-rank tensors $\mu(\vx)$ and $\tilde\e(\vx)$ represent the effective
magnetic permeability and complex permittivity of the corresponding medium;
the latter is $\tilde\e(\vx)=\e(\vx)+i \sigma(\vx)/\omega$, where $\e(\vx)$
and $\sigma(\vx)$ are the (second-rank tensorial) dielectric permittivity
and conductivity.  We assume that $(\vE, \vB)$, $(\vJa, \vMa)$ and
($\tilde\e, \mu)$ in~\eqref{eq:timeharmonicmaxwell} are $\vx$ dependent
with some (weak) regularity of the fields to ensure unique solvability, as
discussed in Section~\ref{sec:variational}.

Equations~\eqref{eq:timeharmonicmaxwell}, interpreted in the strong sense,
hold in appropriate unbounded regions of the $n$-dimensional Euclidean
space, $\mathbb{R}^n$ ($n=2,3$), excluding the set of points comprising the
conducting sheet. We now turn our attention to the requisite boundary
conditions along the sheet. This is modeled as an idealized, oriented
hypersurface $\Sigma$, $\Sigma\subset \mathbb{R}^n$, with unit normal $\vec
\nu$ and effective surface conductivity
$\sigma^\Sigma(\vx)$~\cite{bludov13, hanson08, hanson11}. This
consideration amounts to a jump condition in the tangential component of
the magnetic field while the tangential electric field is continuous,
viz.,~\cite{bludov13}
\begin{align}
  \begin{cases}
    \begin{aligned}
      \vec\nu\times\big\{\big(\mu^{-1}\vB\big)^+-\big(\mu^{-1}\vB\big)^-\big\}
      \Big|_{\Sigma}
      \;&=\;
      \sigma^\Sigma(\vx)\,\big\{(\vec \nu \times\vE)\times\vec \nu\big\}
      \Big|_{\Sigma}~,
      \\[0.2em]
      \vec \nu\times\big\{\vE^+-\vE^-\big\}\Big|_{\Sigma}
      \;&=\; 0~,
    \end{aligned}
  \end{cases}
  \label{eq:jumpcondition}
\end{align}
where $\vMF^{\pm}$ ($\vMF= \vE, \vB$) is the restriction of the
vector-valued solution to either side ($\pm$) of the hypersurface. The
surface conductivity, $\sigma^{\Sigma}(\vx)$, is in principle a second-rank
tensor and is responsible for the creation of the SPP under the appropriate
source $(\vJa, \vMa)$~\cite{bludov13}; see section~\ref{subsec:SPP-intro}.
At the terahertz frequency range in doped graphene, for example, it is
possible that the jump in the tangential component of the magnetic field is
small compared to the magnitude of the field itself~\cite{bludov13}. For
the appropriate polarization and imaginary part of $\sigma^\Sigma$, this
feature may yield a surface wave, the SPP, with a wavelength of the order
of a few microns, much smaller than the free-space
wavelength~\cite{bludov13, hanson08, margetis15}.

In addition, the electromagnetic field $(\vE, \vB)$ must satisfy the
Silver-Müller radiation condition, an extension of the Sommerfeld radiation
condition, if the ambient (unbounded) medium is isotropic~\cite{muller69}.
This amounts to the requirement that $\vMF$ $(\vMF= \vE, \vB)$ approach a
spherical wave uniformly in the radial direction as $|\vx|\to
\infty$ for points at infinity and away from the conducting sheet.
We need to impose
\begin{equation}\label{eq:S-M_cond}
\lim_{|\vx|\to\infty}\{\vB\times \vx-c^{-1}|\vx|\,\vE\}=0~,\
\lim_{|\vx|\to\infty}\{\vE\times \vx+ c|\vx|\, \vB\}=0~,\ x\notin \Sigma~;
\end{equation}
$c$ is the speed of light in the respective medium. In the formulation of
our numerical scheme, we avoid making explicit use of
condition~\eqref{eq:S-M_cond} by using appropriate boundary conditions
together with a PML, which eliminates reflection from infinity.


\subsection{Rescaling}
\label{sec:rescaling}

Loosely following Colton and Kress~\cite{colton83}, as well as
Monk~\cite{monk03}, we introduce a rescaling for the time-harmonic Maxwell
equations (\ref{eq:timeharmonicmaxwell}). The key differences of our
formulation from the above treatments~\cite{colton83,monk03} are:
\begin{itemize}
  \item The additional rescaling of every length scale in our problem by
    the free-space wavelength $2\pi k_0^{-1}:=2\pi
    (\omega\sqrt{\e_0\mu_0})^{-1}$, where $\epsilon_0$ and $\mu_0$ denote
    the {\em vacuum} dielectric permittivity and magnetic permeability,
    respectively. This rescaling recognizes that the typical length scale
    of the SPP is one to two orders of magnitude smaller than the
    corresponding free-space wavelength~\cite{bludov13}; consequently,
    $1/k_0$ is the appropriate macroscopic length scale.
  \item The rescaling of $\vE$, $\vB$, $\vJa$, and $\vMa$ by a
    \emph{typical electric current strength}, $J_0$. In our case, $J_0$ is
    the strength of the prescribed dipole source at location $\vec a$ in
    the $\vec e_i$ direction in Cartesian coordinates:
    \begin{align*}
      \vJa=J_0\,\vec e_i\,\delta(\vx -\vec a).
    \end{align*}
\end{itemize}

Accordingly, we rescale $\mu$ and $\tilde \e$  by $\mu_0$ and $\e_0$,
respectively; cf.~\eqref{eq:timeharmonicmaxwell}:
\begin{alignat}{5}
    \label{eq:rescaling}
    \mu
    &\quad\longrightarrow\quad&
    \mur &=\frac1{\mu_0}\mu,
    &\qquad&
    {\tilde\e}
    &\quad\longrightarrow\quad&
    \er &=&\;
    \frac1{\e_0}\tilde\e.
    \intertext{Furthermore, by use of the free-space wave number,
      $k_0=\omega\sqrt{\e_0\mu_0},$ and the dipole strength, $J_0$, the
      rescaling of the vector fields and coordinates is carried out:}
    \vx
    &\quad\longrightarrow\quad&
    \hat\vx &= k_0\,\vx,
    &\qquad&
    \nabla
    &\quad\longrightarrow\quad&
    \hat\nabla &=&\; \frac1{k_0}\,\nabla,
    \\[0.1em]
    \vJa
    &\quad\longrightarrow\quad&
    \hvJa &= \frac1{J_0}\,\vJa,
    &\qquad&
    \vMa
    &\quad\longrightarrow\quad&
    \hvMa &=&\; \frac {k_0}{\omega\mu_0\,J_0}\,\vMa,
    \\[0.1em]
    \vE
    &\quad\longrightarrow\quad&
    \hvE &= \frac{k_0^2}{\omega\mu_0\,J_0}\,\vE,
    &\qquad&
    \vB
    &\quad\longrightarrow\quad&
    \hvB &=&\; \frac {k_0}{J_0}\,
    \mu^{-1}
    \vB.
\end{alignat}
In addition, the interface conductivity is rescaled as follows:
\begin{align}
    \sigma^\Sigma
    \quad\longrightarrow\quad
    \ssr
    = \sqrt{\frac{\mu_0}{\e_0}}\,\sigma^\Sigma.
\end{align}
Finally, rescaling time-harmonic Maxwell's equations
(\ref{eq:timeharmonicmaxwell}) results in the following system:
\begin{align}
  \begin{cases}
    \begin{aligned}
      -i\mur\hvB+\hat\nabla\times\hvE \;&=\; -\hvMa,
      \\[0.1em]
      \hat\nabla\cdot\big(\mur\hvB\big) \;&=\; \frac 1{i}
      \hat\nabla\cdot\hvMa,
      \\[0.1em]
      -i\er\hvE-\hat\nabla\times\hvB \;&=\; -\hvJa.
      \\[0.1em]
      \hat\nabla\cdot\big(\er\hvE\big) \;&=\;\frac1i \nabla\cdot\hvJa.
    \end{aligned}
  \end{cases}
  \label{eq:timeharmonicmaxwellrescaled}
\end{align}
To lighten the notation, the hat ($\,\hat\ \,$) on top of a rescaled
quantity will be omitted in the remainder of this paper. (This
simplification avoids confusion of the rescaled quantities with the Fourier
transforms of fields invoked in Section~\ref{sec:analytical}).


\subsection{Variational statement}
\label{subsec:variational}

\begin{figure}[t]
  \centering
  \includegraphics{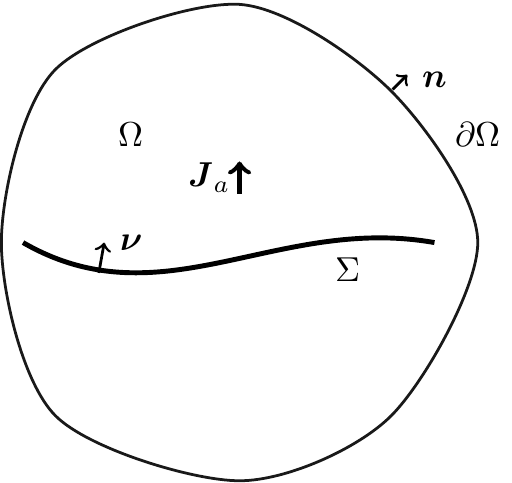}
  \caption{
  Schematic of the computational domain, $\Omega$, with boundary
  $\partial\Omega$ and outer normal $\vec n$. An electric Hertzian dipole,
  $\vJa$, is situated above a prescribed hypersurface, $\Sigma$.
  }
  \label{fig:domain}
\end{figure}

Let $\Omega\subset\R^n$ ($n=2,3$) be a simply connected and bounded domain
with Lipschitz-continuous and piecewise smooth boundary, $\partial\Omega$.
Further, let $\Sigma$ be an oriented, Lipschitz-continuous, piecewise
smooth hypersurface. Fix a normal field $\vn$ on $\Sigma$ and let $\vec n$
denote the outer normal vector on $\partial\Omega$ ; see
Figure~\ref{fig:domain}. Substituting $\vB$ from the first equation of
(\ref{eq:timeharmonicmaxwellrescaled}) into the third equation yields
\begin{align}
  \nabla\times\big(\muri\nabla\times\vE\big)
  -\er\vE
  \;=\; i\,\vJa - \nabla\times\big(\muri\vMa\big).
  \label{eq:2ndorderequation}
\end{align}
Multiplying (\ref{eq:2ndorderequation}) by the complex conjugate,
$\bar\vp$, of a smooth test function $\vp$ and integrating by parts in
$\Omega\setminus\Sigma$ lead to
\begin{multline}
  \int_\Omega (\muri\nabla\times\vE)\cdot(\nabla\times\bar\vp)\dx
  \;\;-\;
  \int_\Omega \er\vE\cdot\bar\vp\dx
  \\
  -\int_\Sigma\big[\vn\times(\muri\nabla\times\vE+\muri\vMa)\big]_\Sigma
  \cdot\bar\vp_T\dox
  \\
  +\;
  \int_{\partial\Omega}\big(\vn\times(\muri\nabla\times\vE+\muri\vMa)\big)
  \cdot\bar\vp_T\dox
  \\
  =\;\;i\int_\Omega\vJa\cdot\bar\vp\dx -
  \int_\Omega\muri\vMa\cdot\nabla\times\bar\vp\dx,
  \label{eq:integratedbyparts}
\end{multline}
where the subscript $T$ above denotes the tangential part of the respective
vector, $\vMF_T = (\vn\times\vMF)\times\vn$, and $[\,.\,]_\Sigma$ denotes
the jump over $\Sigma$ with respect to $\vn$, viz.,
\begin{align}
  \big[\vMF\big]_\Sigma(\vx)
  :=\lim_{s\searrow0}\big(\vMF(\vx+s\vn)-\vMF(\vx-s\vn)\big)\qquad (\vx\in\Sigma).
\end{align}
For the computational domain $\Omega$, an \emph{absorbing boundary
condition} at $\partial\Omega$ is imposed:
\begin{align}
  \vn\times\vB+\sqrt{\muri\er}\vE_T = 0\qquad (\vx\in\partial\Omega).
  \label{eq:absorbingbc}
\end{align}
The last boundary condition is obtained by using a first-order
approximation of the Silver-Müller radiation condition,
equation~\eqref{eq:S-M_cond}, truncated at $\partial\Omega$~\cite{jin08}.
We will occasionally refer to \eqref{eq:absorbingbc} as the
\emph{first-order} absorbing boundary condition. We assume that $\muri$ and
$\er$ are scalar functions such that the square root in
(\ref{eq:absorbingbc}) is well defined. In our numerical computation, we
combine the above absorbing boundary condition with a PML; see
Section~\ref{sub:pollutionandpml}.

An advantage of variational formulation (\ref{eq:integratedbyparts}) is
that the jump condition over the conducting sheet can be expressed as a
\emph{weak discontinuity}. Rewriting jump condition
(\ref{eq:jumpcondition}) as well as absorbing boundary condition
(\ref{eq:absorbingbc}) in terms of $\vB$ and $\vMa$ by utilizing the first
equation of (\ref{eq:timeharmonicmaxwellrescaled}) yields
\begin{align}
  \big[\vn\times(\muri\nabla\times\vE+\muri\vMa)\big]_\Sigma &=
  i\,\ssr\vE_T\quad\text{on }\Sigma, \label{eq:jumpconditionrescaled}
  \\[0.2em]
  \vn\times(\muri\nabla\times\vE+\muri\vMa) &= -i\sqrt{\muri\er}\,\vE_T
  \quad\text{on }\partial\Omega.
\end{align}
The last two relations allow us to enforce the jump and boundary conditions
weakly by simply substituting them into (\ref{eq:integratedbyparts}). We
summarize our main result below.
\medskip

\begin{proposition}[Variational formulation]
  \label{prop:existence}

  In order to ensure unique solvability, we assume that $\ssr\in
  L^\infty(\Sigma)^{3\times 3}$ is matrix-valued and symmetric, with
  semi-definite real and complex part. Further, let $\er$ be a smooth
  scalar function and $\muri$ be a constant scalar such that
  \begin{itemize}
    \item
      $\text{Im}\,\big(\er)=0$, or $\text{Im}\,\big(\er)\ge c>0$ in
      $\Omega$,
    \item
      $\sqrt{\muri\er}$ is real-valued and strictly positive on
      $\partial\Omega$.
  \end{itemize}
  Define a Hilbert space (cf.\ \cite[Th.\,4.1]{monk03})
  \begin{align*}
    \vec X(\Omega)=\Big\{\vp\in \vec H(\text{curl};\Omega)\;:\;
    \vp_T\big|_\Sigma\in L^2(\Sigma)^3,\;
    \vp_T\big|_{\partial\Omega}\in L^2(\partial\Omega)^3
    \Big\}
  \end{align*}
  equipped with the norm $\|\vp\|^2_{\vec X}=\|\nabla\times\vp\|^2+
  \|\vp_T\|_{L^2(\Sigma)}^2 + \|\vp_T\|_{L^2(\Omega)}^2$. In the above,
  $\vec H(\text{curl}; \Omega)$ denotes the space of vector-valued,
  measurable and square integrable functions whose (distributive) curl
  admits a representation by a square integrable function.

  The rescaled, weak formulation of (\ref{eq:timeharmonicmaxwell}) with
  jump condition (\ref{eq:jumpcondition}) and absorbing boundary condition
  (\ref{eq:absorbingbc}) can be stated as follows:
  Find $\vE\in\vec X(\Omega)$,
  such that
  \begin{align}
    A(\vE,\vp) = F(\vp),
    \label{eq:variationalformulation}
  \end{align}
  for all $\vp\in\vec X(\Omega)$. It admits a unique solution. The
  sesquilinear form and the right-hand side are given by
  \begin{align}
    A(\vE,\vp) &:=
    \int_\Omega (\muri\nabla\times\vE)\cdot(\nabla\times\bar\vp)\dx
    \;\;-\;
    \int_\Omega\er\vE\cdot\bar\vp\dx
    \\\notag
    &\qquad
    -\, i\,
    \int_{\Sigma}(\ssr\vE_T)\cdot\bar\vp_T\dox
    \;\;-\; i\,
    \int_{\partial\Omega}\sqrt{\muri\er}\vE_T\cdot\bar\vp_T\dox,
    \\[0.1em]
    F(\vp)
    &:=\;\;
    i\int_\Omega\vJa\cdot\bar\vp\dx -
    \int_\Omega\muri\vMa\cdot\nabla\times\bar\vp\dx.
    \label{eq:variationalformulation2}
  \end{align}
\end{proposition}
\medskip

\begin{proof}%
  The existence result for time-harmonic Maxwells equations with an
  absorbing boundary condition~\cite{colton83,colton98,monk03} can be
  applied almost directly to problem (\ref{eq:variationalformulation}); the
  additional interface integral in $A(\vE,\vp)$,
  \begin{align}
    -\,i\,\int_\Sigma\ssr\vE_T\cdot\bar\vp_T\dox.
  \end{align}
  requires a careful discussion. For this we split the integral into two
  contributions,
  \begin{align}\label{eq:surface-int}
    -\,i\,\int_\Sigma\ssr\vE_T\cdot\bar\vp_T\dox
    =
    -\,i\,\int_\Sigma\text{Re}\,\big(\ssr\big)\vE_T\cdot\bar\vp_T\dox
    +
    \int_\Sigma\text{Im}\,\big(\ssr\big)\vE_T\cdot\bar\vp_T\dox.
  \end{align}
  The first term on the right-hand side of~\eqref{eq:surface-int} can be
  treated similarly to the absorbing boundary condition on $\partial\Omega$
  (cf. \cite[Sec.\,4.5]{monk03}). For the second term
  in~\eqref{eq:surface-int}, involving $\text{Im}\,\big(\ssr\big)$, it
  holds that
  \begin{align}
    \int_\Sigma\text{Im}\,\big(\ssr\big)\vE_T\cdot\bar\vE_T\dox
    \ge 0.
  \end{align}
  Thus, this term does not affect any proof based on showing coercivity of
  a modified sesquilinear form and using Fredholm's alternative on a
  compact perturbation.

  In order to prove uniqueness we follow the strategy
  in~\cite[Sec.\,4.6]{monk03}. First, test the homogeneous equation
  $A(\ve,\vp) = 0$ with $\ve$ itself and take
  the imaginary part,
  \begin{align*}
    \int_\Omega(\text{Im}(\er)\,\ve\big)\cdot\bar\ve\dx
    +
    \int_{\Sigma}(\text{Im} (\ssr)\,\ve_T)\cdot\bar\ve_T\dox
    +
    \int_{\partial\Omega}\sqrt{\muri\er}\ve_T\cdot\bar\ve_T\dox = 0.
  \end{align*}
  This immediately implies $\ve_T=0$ on $\partial\Omega$ and $\Sigma$. The
  (nontrivial) conclusion $\ve=0$ now follows verbatim by the proof
  of~\cite[Th.\,4.12]{monk03}.
\end{proof}
\medskip

The following remarks are in order.

\begin{remark}[2D model of a conducting sheet]
  \label{remark:2dmodel}
  The 3D problem description given by (\ref{eq:variationalformulation})
  readily translates into a corresponding problem in 2D: Assume that the
  interface $\Sigma$, the electric field $\vE$, the permeability $\mur$, as
  well as the permittivity $\er$ and surface conductivity $\ssr$ are
  translation- and mirror-invariant in the $z$-direction. Thus, the term
  $\nabla\times\vE$ and, consequently, the $\vB$ field have nonzero
  components only in the $z$-direction. Hence, we can rewrite
  (\ref{eq:variationalformulation}) with the 2D curl
  \begin{align}
    \begin{aligned}
      \nabla\times\vMF := \partial_x\mathcal{F}_y - \partial_y\mathcal{F}_x.
    \end{aligned}
  \end{align}
\end{remark}
\medskip

\begin{remark}[TM polarization in 2D model]
  In the 2D version of Maxwell's equations with a vertical electric dipole
  (see Section~\ref{sec:analytical}), the magnetic field, $\vB$, given by
  \begin{align}
    \vB=\muri\big(\vMa+\nabla\times\vE\big),\qquad \vMa\equiv 0~,
  \end{align}
  only has a $z$-component when viewed as a vector field in $\mathbb{R}^3$.
  Thus, this field is parallel to the hypersurface $\Sigma$ and orthogonal
  to the computational domain, $\Omega$ (which is part of the $xy$-plane).
  Consequently, the SPP in this 2D setting, and in the corresponding
  numerical simulation, has the desired TM polarization~\cite{bludov13}.
\end{remark}


\section{Exact solution for 2D infinite conducting sheet}
\label{sec:analytical}
\begin{figure}
  \begin{center}
  \includegraphics{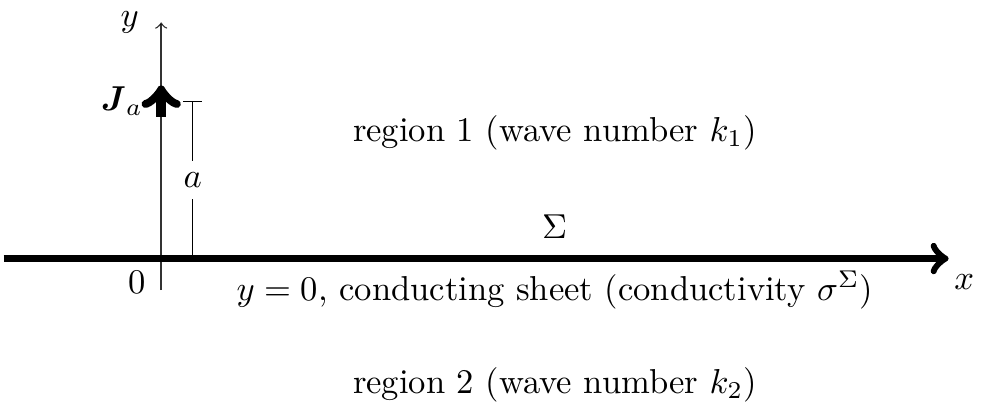}
  \end{center}

  \caption{Schematic of a vertical electric dipole at distance $a$ from
    conducting sheet $\Sigma$ in 2D. The dipole has current density
    $\vJa=J_0 \delta(\vx-\vec a)\,e_y$, where  $\vx= x e_x+ y e_y$ and
    $\vec a= a\,e_y$. The sheet lies in $y=0$, separates the space into
    region $1$ (half space $\{y>0\}$) with wave number $k_1$ and region $2$
    ($\{y<0\}$) with wave number $k_2$, and has surface electric
    conductivity $\sigma^\Sigma$.}
 \label{fig:geometry}
\end{figure}

In this section, we derive an exact solution to the three-dimensional (3D)
version of boundary value
problem~\eqref{eq:timeharmonicmaxwell}--\eqref{eq:S-M_cond} in the case
with a 2D vertical electric dipole radiating at frequency $\omega$ over an
infinite conducting sheet embedded in an isotropic and homogeneous space of
(scalar) magnetic permeability $\mu$; see Figure~\ref{fig:geometry}. For
physical clarity, we invoke the vector-valued
electromagnetic field without rescaling, unless we state otherwise.

The sheet separates the space into region 1 ($\{y>0\}$) with wave number
$k_1$ and region 2 ($\{y<0\}$) with wave number $k_2$; $k_j^2=\omega^2
\tilde\epsilon_j\mu$ where $\tilde\epsilon_j$ is the complex permittivity
($j=1,\,2$). We assume that the surface conductivity, $\sigma^\Sigma$, of
the sheet is a scalar constant. Note that we assign different complex
permittivities to each half-space (regions 1, 2). In the end of this
section, we set them equal.

The dipole has current density $\vJa= J_0\,\delta(\vx-\vec a)\,e_y$ where
$\vx=x\,e_x+y\,e_y$ and $\vec a= a\,e_y$; $e_s$ is a unit Cartesian vector
($s= x, y$). Now define the Fourier transform, $\widehat{\vF}_j(\xi, y)$,
of the vector-valued field $\vF_j(x, y)$ ($\vF=\vE,\,\vB$) in region $j$
through the integral formula
\begin{equation}
  \label{eq:EB-FT}
  \vF_j(x,y)=\frac{1}{2\pi}\int_{\mathbb{R}}d\xi\ \widehat{\vF}_j(\xi, y)\,
  e^{i \xi x}~.
\end{equation}
The transformation of 3D Maxwell's equations~\eqref{eq:timeharmonicmaxwell}
in region $j$ yields
\begin{align}
  & -i\xi \widehat{E}_{jy}+\frac{\partial \widehat{E}_{jx}}{\partial y}
    =-i\omega \widehat{B}_{jz}~,\label{eq:maxwell_1}\\
  & -\frac{\partial\widehat{B}_{jz}}{\partial y}
    =\frac{ik_j^2}{\omega}\widehat{E}_{jx}~,\ -i\xi\widehat{B}_{jz}
    =-\frac{ik_j^2}{\omega}\widehat{E}_{jy}+\mu\, \delta (y-a) ~,
    \label{eq:maxwell_2}
\end{align}
where we set $E_{jz}\equiv 0$, $B_{jx}\equiv 0,$ and $B_{jy}\equiv 0$ by
symmetry; $F_{jz}$ ($F= E,\,B$) is the component perpendicular to the
$xy$-plane.

Equations~\eqref{eq:maxwell_1} and~\eqref{eq:maxwell_2} furnish the
differential equation
\begin{equation*}
  \frac{\partial^2\widehat{B}_{jz}}{\partial y^2}+\beta_j^2\,\widehat{B}_{jz}
  =-i\xi \mu\,\delta(y-a)~,\quad \beta_j^2:=k_j^2-\xi^2~,
\end{equation*}
which has solution
\begin{equation}\label{eq:soln-Bz}
  \widehat{B}_{jz}(\xi,y)=\left\{\begin{array}{lr}
  \displaystyle{C_{>}\,e^{i\beta_1 y}-\frac{\xi\mu}{2\beta_1}\,
  e^{i\beta_1|y-a|}}~,&\ \mbox{if}\ j=1\ (y>0)~;\\
  C_{<}\, e^{-i\beta_2 y}~,&\ \mbox{if}\ j=2\ (y<0)~.\\
  \end{array}\right.
\end{equation}
This solution is consistent with radiation condition~\eqref{eq:S-M_cond}
provided
\begin{equation}\label{eq:beta-cond}
\text{Im}\beta_j(\xi)>0\quad (j=1,\,2)~.
\end{equation}

Furthermore, we apply conditions~\eqref{eq:jumpcondition} in order to
determine integration constants $C_>$ and $C_<$. Specifically, we impose
\begin{equation}
  \big(\widehat{B}_{1z}-\widehat{B}_{2z}\big)\bigl|_{y=0}=\mu\,\sigma^\Sigma
  \widehat{E}_{1x}\bigl|_{y=0}~,\quad
  \big(\widehat{E}_{1x}-\widehat{E}_{2x}\big)\bigl|_{y=0}=0~.
\end{equation}
Accordingly, by~\eqref{eq:maxwell_2} and~\eqref{eq:soln-Bz} we find
\begin{equation*}
  C_>=-\frac{\xi\mu}{2\beta_1}\,\frac{k_2^2\beta_1-k_1^2\beta_2+
  \omega\mu\sigma^\Sigma \beta_1\beta_2}{k_2^2\beta_1+k_1^2\beta_2+
  \omega\mu\sigma^\Sigma \beta_1\beta_2 }e^{i\beta_1 a}~,\
  C_<=-\frac{\mu k_2^2\xi\,e^{i\beta_1a}}{k_2^2\beta_1+k_1^2\beta_2+
  \omega\mu\sigma^\Sigma\beta_1\beta_2}~.
\end{equation*}

We can now write down all field components in view
of~\eqref{eq:EB-FT}--\eqref{eq:soln-Bz}. In particular:
\begin{align}
  E_{1x}(x,y)&=\frac{\omega\mu}{4\pi k_1^2}\int_{-\infty}^\infty d\xi\
  \xi\Biggl[\frac{k_2^2\beta_1-k_1^2\beta_2+\omega\mu\sigma^\Sigma
  \beta_1\beta_2}{k_2^2\beta_1+k_1^2\beta_2+\omega\mu\sigma^\Sigma
  \beta_1\beta_2 } e^{i\beta_1(y+a)}
  \label{eq:E1x-ex}\\\nonumber
  &\qquad +\text{sgn}(y-a)\, e^{i\beta_1 |y-a|}\Biggr]e^{i\xi x}
  \quad (y>0)~;
  \\
  E_{2x}(x,y)&=-\frac{\omega\mu}{2\pi}\int_{-\infty}^\infty d\xi\
  \xi\beta_2\,\frac{e^{i\beta_1a-i\beta_2y}e^{i\xi x}}
  {k_2^2\beta_1+k_1^2\beta_2+\omega\mu\sigma^\Sigma\beta_1\beta_2}\quad (y<0)~.
  \label{eq:Ex-ex}
\end{align}
In the above, $\text{sgn}(y)=1$ if $y>0$ and $\text{sgn}(y)=-1$ if $y<0$.

Note that all field components have
Fourier transforms defined in the $\xi$-plane with: (i) branch points at
$\xi=\pm k_j$ ($j=1, 2$); and (ii) simple poles at points where
\begin{equation}
  \label{eq:SPP-dispersion}
  k_2^2\beta_1(\xi)+k_1^2\beta_2(\xi)+
  \omega\mu\sigma^\Sigma\beta_1(\xi)\beta_2(\xi)=0~,
\end{equation}
under condition~\eqref{eq:beta-cond} which defines the appropriate branch
of the multiple-valued $\beta_j(\xi)$. Equation~\eqref{eq:SPP-dispersion}
is the dispersion relation for the SPP~\cite{raether86,margetis15}. In
particular, if $k_1=k_2=k$ and $|\omega\mu\sigma^\Sigma|\ll |k|$,
\eqref{eq:SPP-dispersion} reduces to~\eqref{eq:simple_disp-reln} with
$\xi=k_x$.

This observation motivates the following definition within the 2D model~\cite{margetis15}.
\medskip

\begin{definition}[SPP in the 2D setting]
  For an infinite conducting sheet, the SPP is
  identified with the part of the electromagnetic field equal to the
  contribution to the Fourier integrals of the simple pole $\xi=k_m$,
  $\text{Im} k_m>0$, that solves~\eqref{eq:SPP-dispersion}
  under~\eqref{eq:beta-cond}.\label{def:SPP}
\end{definition}

For the sake of completeness, we conclude this section by focusing on the
case with $k_1=k_2=k$ (identical half-spaces). In particular, we provide
explicit expressions for two distinct physical contributions to the
$x$-component of the electric field on the sheet ($y=0$), which is
continuous across the sheet. After suitable rescaling of the variables and
parameters (Section~\ref{sec:rescaling}), by~\eqref{eq:Ex-ex} the pole
contribution to $E_{x}:=E_{1x}(x,0)=E_{2x}(x,0)$ takes the form
\begin{align}
  \label{eq:polecontributionrescaled}
  E_{x,\text{p}}
  = -2i\frac{\mur\er}{(\ssr)^2}\exp\big[i\kmr x - (2i/\ssr)\,a_r]~,
\end{align}
where $\kmr=k_m/k_0$ and $a_r=k_0 a$; cf.~\eqref{eq:SPP-dispersion}.
Equation~\eqref{eq:polecontributionrescaled} expresses the SPP pertaining
to the tangential electric field (see Definition~\ref{def:SPP}).

A separate physical contribution expresses the radiation part of the
scattered field, $\vE^{\rm sc}=\vE-\vE^{\rm pr}$, where $\vE^{\rm pr}$ is
the (primary) dipole field in the absence of the sheet. To single out
this contribution for the $x$-component, $E_x$, we choose to integrate in
the $\xi$ (Fourier) variable around the branch cut emanating from the
branch point $\xi=k$ after removal of the primary field component. The
result reads
\begin{multline}
  \label{eq:branchcutcontributionrescaled}
  E_{x,\text{bc}}^{\text{sc}} =
  \frac1{4\pi}\frac{1}{\ssr}
  \Bigg\{
    \int_0^1\text{d}\xi\ \xi\sqrt{1-\xi^2}\,e^{i\sqrt{\mur\er}x\,\xi}
    \,\frac1{\xi^2+4\frac{\mur\er}{(\ssr)^2}-1}
    \\
    \cdot\Big(4\sqrt{\mur\er}\cos\big(\sqrt{\mur\er}a\sqrt{1-\xi^2}\big)-
         2i\ssr\sqrt{1-\xi^2}\sin\big(\sqrt{\mur\er}a\sqrt{1-\xi^2}\big)\Big)
    \\
    -
    \int_0^\infty\text{d}\varsigma\ \varsigma\sqrt{1+\varsigma^2}\,e^{i\sqrt{\mur\er}x\,\varsigma}
    \,\frac1{\varsigma^2-4\frac{\mur\er}{(\ssr)^2}-1}
    \\
    \cdot\Big(4\sqrt{\mur\er}\cos\big(\sqrt{\mur\er}a\sqrt{1+\varsigma^2}\big)-
         2i\ssr\sqrt{1+\varsigma^2}\sin\big(\sqrt{\mur\er}a\sqrt{1+\varsigma^2}\big)\Big)
  \Bigg\}.
\end{multline}
\begin{figure}[!t]
  \centering
  \subfloat[]{
    \includegraphics{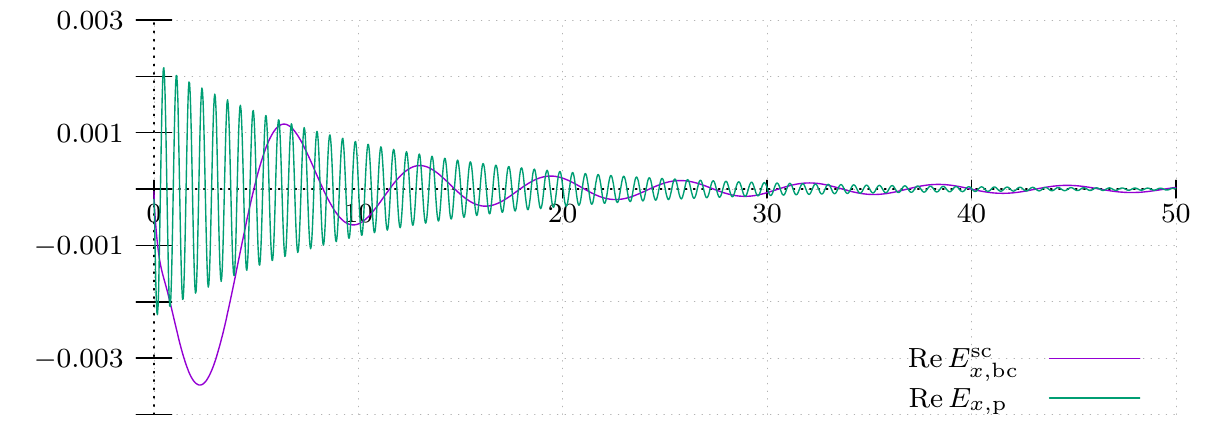}
  }

  \subfloat[]{
    \includegraphics{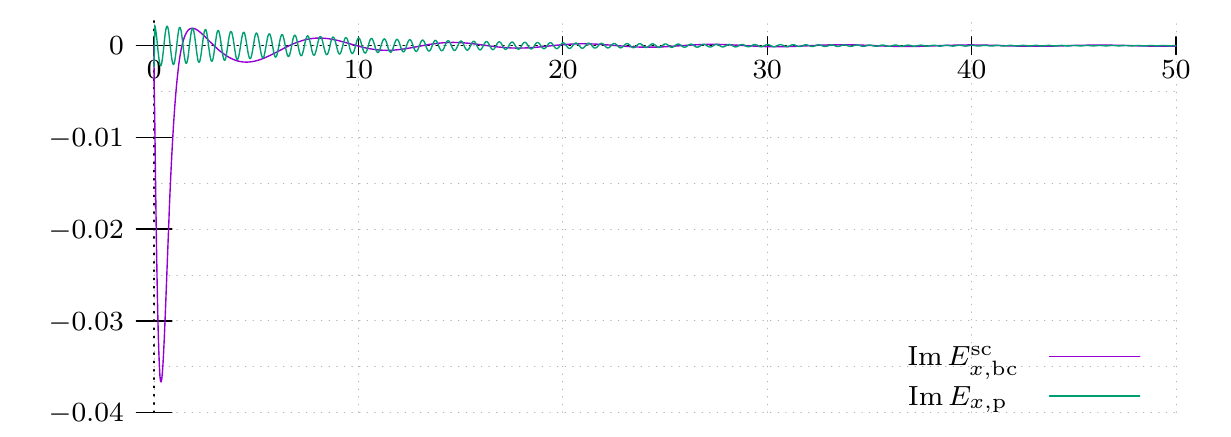}
  }
  \caption{
    Real (a) and imaginary (b) part of branch-cut contribution
    \eqref{eq:branchcutcontributionrescaled} and SPP
    \eqref{eq:polecontributionrescaled} for a dipole with (rescaled)
    elevation distance $a=1.0$ and SPP wave number $\kmr=10.0+0.1i$. The
    ambient space has $\er=1$ and $\mur=1$.}
  \label{fig:pole-vs-branchcut}
\end{figure}

Figure~\ref{fig:pole-vs-branchcut} shows a comparison of the two
contributions, $ E_{x,\text{p}}$ and $E_{x,\text{bc}}^{\text{sc}}$, for a
typical choice of parameters ($a=1.0$,
$\sigma^\Sigma=2.00\times10^{-3}+0.200i$) by use of
formulas~\eqref{eq:polecontributionrescaled}
and~\eqref{eq:branchcutcontributionrescaled} with $\mur=1$ and $\er=1$. We
observe that the SSP (pole contribution) dominates in the (rescaled) range
$10\le x\le 25$ but is eventually dominated by the branch-cut contribution
which has the (slower) algebraic decay.


\section{Numerics: Discretization scheme}
\label{sec:numerics}

In this section, variational formulation (\ref{eq:variationalformulation})
is discretized on a non-uniform quadrilateral mesh with higher-order,
curl-conforming Nédélec elements~\cite{nedelec86,nedelec01,bokil15}.
Such a conforming discretization is an ideal choice for the problem at
hand. The interface with the (weak) discontinuity can be aligned with the
mesh and away from it; and the regularity of the solution leads to high
convergence rates. Key ingredient of our treatment is the use of an
appropriately defined PML (Section \ref{sec:pml}), as well as local mesh
refinement based on a posteriori error estimates (Section
\ref{sec:aposteriori}).

Let $X_h(\Omega)\subset X(\Omega)$ be a finite element subspace spanned by
Nédélec elements. We will in particular use second-order Nédélec elements
in the numerical computations. Then, under the assumption of a sufficiently refined
initial mesh, the variational equation to find $E_h\in X_h(\Omega)$ such that
\begin{align}
  A(\vE_h,\vp) = F(\vp),
  \quad\forall\vp\in X_h(\Omega),
\end{align}
is uniquely solvable \cite[Section 7.2]{monk03}.
\begin{remark}
  From an approximation theory point of view, we expect the convergence
  \begin{align}
    \big\|\vE-\vE_h\big\|_X\sim\mathcal{O}(h^2)=\mathcal{O}(\#\text{dofs})
    \label{eq:convergenceorder}
  \end{align}
  for second-order Nédélec elements and under the conditions of
  Proposition~\ref{prop:existence}, i.\,e., quadratic convergence in terms
  of a uniform refinement parameter $h$, or linear convergence in the number of
  degrees of freedom. This is evidenced by our numerical results presented
  in Section~\ref{sec:computationalresults}. We refrain from stating a
  rigorous convergence result at this point because the use of non-uniform,
  locally refined meshes with an additional approximation of a curved
  boundary significantly complicates the classical approximation theory
  (see \cite[Section 7.3]{monk03} and references therein).
\end{remark}

The interface jump condition on the sheet introduces a pronounced two-scale
character to the problem that needs special numerical treatment. Depending
on the surface electric conductivity, $\sigma^\Sigma$, of the sheet, the
observed SPP may have a rescaled wave number, $\kmr$, that is one to two
orders of magnitude higher than that produced by the dipole radiation in
free space (which has rescaled wave number $k_r=1$). This fact has
important consequences with respect to the refinement strategy and boundary
conditions. In particular:
\begin{itemize}
  \item
    First-order absorbing boundary conditions alone are in principle not
    well suited for conducting sheets sustaining SPPs. These boundary
    conditions lead to a significant suppression of the SPP amplitude. This
    is especially an issue for configurations where the SPP given by
    \eqref{eq:polecontributionrescaled} has a significantly smaller
    amplitude than the branch-cut contribution
    \eqref{eq:branchcutcontributionrescaled}.
  \item
    A much finer minimal mesh refinement is necessary near the interface
    $\Sigma$ in order to resolve the highly oscillatory SPP. In addition,
    failure to sufficiently resolve the interface in the \emph{whole}
    computational domain results in a suppression of the SPP due to
    non-local cancellation effects (\emph{pollution effect}); see
    Section~\ref{sub:pollutionandpml} for computational examples.
\end{itemize}
One of the major advantages of a finite element approach for discretizing
variational formulation (\ref{eq:variationalformulation}) is the fact that
no regularization of the interface by a layer with artificial thickness is
necessary. Instead, the sheet can be \emph{directly approximated} as a
lower-dimensional interface.

In the remainder of this section, we discuss our choice of a PML to remedy
the negative effect of the absorbing boundary condition on the SPP
amplitude. Further, a strategy for adaptivity and \emph{local mesh
refinement} is introduced, which is based on \emph{a~posteriori} error
control combined with a fixed (\emph{a priori}) local refinement of the
interface.


\subsection{Perfectly matched layer}
\label{sec:pml}

In this subsection, we carry out a construction of a PML~\cite{berenger94,
chew94, zhao96} for the rescaled Maxwell equations with a jump condition.
The concept of a PML was pioneered by Bérenger~\cite{berenger94}. It is
essentially a layer with modified material parameters ($\er$, $\mur$)
placed near the boundary such that all outgoing electromagnetic waves decay
exponentially with no ``artificial'' reflection due to truncation of the
domain. The PML is an indispensable tool for truncating unbounded domains
for wave equations and often used in the numerical approximation of
scattering problems~\cite{monk03,bao10,chew94,berenger94,zhao96}.

We use an approach for a PML for time-harmonic Maxwell's equations outlined
by Chew and Weedon~\cite{chew94}, as well as Monk~\cite{monk03}. The idea
is to use a formal change of coordinates from the computational domain
$\Omega\subset\R^3$ with real-valued coordinates
to a domain
$\acute\Omega\subset\{z\in\mathbb{C}\::\:\text{Im\,}z\ge0\}^3$
with complex-valued coordinates and non-negative imaginary
part~\cite{monk03}; and transform back to the real-valued domain. This
results in modified material parameters $\er$, $\mur$ and $\ssr$ within the
PML.

\begin{figure}[t]
  \centering
  \includegraphics{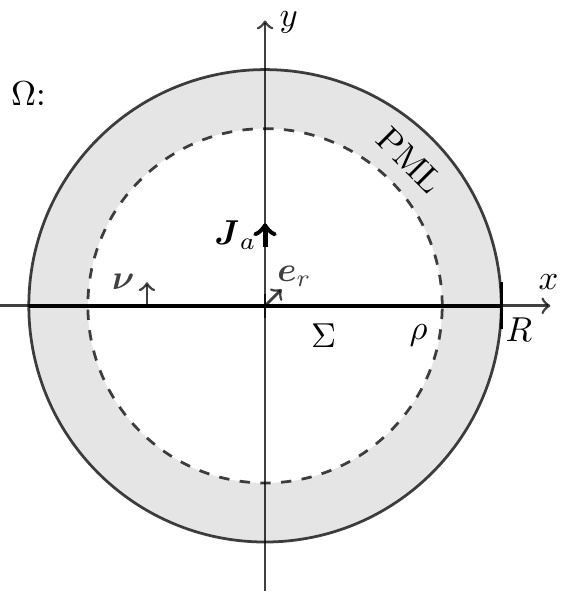}
  \caption{The computational domain $\Omega$ with a vertical dipole $\vJa$
    situated above a planar conducting sheet $\Sigma$ located on the
    $x$-axis.
  }
  \label{fig:graphenesheet}
\end{figure}

We assume that a PML is imposed as a concentric spherical shell in a small
outer region near the boundary $\partial\Omega$; see
Figure~\ref{fig:graphenesheet}. The transformation is chosen to  be a
function of distance $\rho$ in radial ($\vec e_r$-) direction from the
origin. Furthermore, assume that the normal field $\vn$ of $\Sigma$ is
orthogonal to $\vec e_r$ (for $\vx\in \Sigma \cap {\rm (PML)}$) and that
$\vJa\equiv0$ and $\vMa\equiv0$ within the PML. Introduce a change of
coordinates $\Omega\to \acute\Omega$ with
\begin{align}
  \vx \quad\mapsto\quad \acute \vx =
  \begin{cases}
    \begin{aligned}
      &\vx + i \vec e_r \int_\rho^r s(\tau)\,\text{d}\tau
      &\quad &\text{if }r\ge \rho,
      \\[0.1em]
      &\vx & &\text{otherwise},
    \end{aligned}
  \end{cases}
  \label{eq:pmltransformation}
\end{align}
with $r=\vec e_r\cdot\vx$ and an appropriately chosen, nonnegative scaling
function $s(\tau)$~\cite{monk03}. This prescription leads to a modified
system of Maxwell's equations defined on $\acute\Omega$, which takes the
following (rescaled) form within the PML:
\begin{align}
  \begin{cases}
    \begin{aligned}
      \acute\nabla\times\big(\muri\acute\nabla\times\avE\big) -\er\avE
      \;&=\;
      0,
      \\[0.1em]
      \big[\vn\times(\muri\acute\nabla\times\avE)\big]_\Sigma
      \;&=\;
      i\,\ssr\avE_T\quad\text{on }\Sigma,
    \end{aligned}
  \end{cases}
  \label{eq:pmlsystemcomplex}
\end{align}
where the accent on top of  $\vE$ and $\nabla$ indicates the dependence on
as well as the differentiation with respect to $\acute\vx$. Next, we
transform (\ref{eq:pmlsystemcomplex}) back from $\acute\Omega$ to $\Omega$
with the help of diffeomorphism (\ref{eq:pmltransformation}). It follows
that~\cite{monk03}
\begin{align}
  \acute\nabla\times\acute\vMF = A(\nabla\times B\vMF),
  \qquad
  \vn\times\acute\vMF = C(\vn\times B\vMF),
\end{align}
for $\vn$ orthogonal to $\vec e_r$. In the above, we introduce the
$3\times3$ matrices
\begin{gather}
  A=T^{-1}_{\vec e_x\vec e_r}\text{diag}\,
  \Big(\frac{1}{\bar d^2},\frac{1}{d\bar d},\frac{1}{d\bar d}\Big)
  T_{\vec e_x\vec e_r},
  \quad
  B=T^{-1}_{\vec e_x\vec e_r}\text{diag}\,\big(d,\bar d,\bar d\big)
  T_{\vec e_x\vec e_r},
  \\\notag
  C=T^{-1}_{\vec e_x\vec e_r}\text{diag}\,\Big(\frac{1}{\bar d},\frac{1}{\bar d},
  \frac{1}{d}\Big)
  T_{\vec e_x\vec e_r},
\end{gather}
where
\begin{align}
  \label{eq:d}
  d=1+i\,s(r), \quad
  \bar d=1+i/r\int_\rho^r s(\tau)\,\text{d}\tau,
\end{align}
and $T_{\vec e_x\vec e_r}$ is the rotation matrix that rotates $\vec e_r$
onto $\vec e_x$. Thus, applying the rescaling
\begin{align}
  \begin{cases}
    \begin{aligned}
      \muri
      &\quad\longrightarrow\quad
      \breve\mu_r^{-1} = B\muri A,
      \\[0.1em]
      \er
      &\quad\longrightarrow\quad
      \breve\e_r =\; A^{-1}\er B^{-1},
      \\[0.1em]
      \ssr
      &\quad\longrightarrow\quad
      \tissr
      = C^{-1}\ssr B^{-1},
    \end{aligned}
  \end{cases}
  \label{eq:pmlcoefficients}
\end{align}
to (\ref{eq:pmlsystemcomplex}) leads to the system
(with $\breve\vE := B\,\acute\vE$)
\begin{align}
  \begin{cases}
    \begin{aligned}
    \nabla\times\big(\breve\mu_r^{-1}\nabla\times\breve\vE\big)
    -\breve\e_r\breve\vE
      \;&=\;
      0,
      \\[0.2em]
      \big[\vn\times(\breve\mu_r^{-1}\nabla\times\breve\vE)\big]_\Sigma
      \;&=\;
      i\,\tissr\breve\vE_T\quad\text{on }\Sigma,
    \end{aligned}
  \end{cases}
  \label{eq:maxwellpml}
\end{align}
within the \emph{real-valued} domain $\Omega$. Note that outside of the PML
matrices $A$, $B$, and $C$ simply reduce to the unit matrix. Thus,
$\vE=\breve\vE$ and \eqref{eq:maxwellpml} is identical to
\eqref{eq:2ndorderequation} with jump condition
\eqref{eq:jumpconditionrescaled}. The modified equations for the PML can be
implemented by suitably \emph{replacing} $\er$, $\muri$, and $\ssr$ by
their counterparts according to (\ref{eq:pmlcoefficients}) within the PML.
For the 2D model discussed in Remark~\ref{remark:2dmodel}, the scalar
$\muri$ is transformed via
\begin{align}
  \muri
  \quad\longrightarrow\quad
  \breve\mu_r^{-1} = \frac{\muri}d,
\end{align}
where $d$ is given by \eqref{eq:d}.


\subsection{A posteriori error estimation and local refinement}
\label{sec:aposteriori}
One of the computationally challenging aspects of the numerical simulation
is the \emph{two-scale} behavior of problem
\eqref{eq:variationalformulation}. A much finer minimal mesh refinement is
necessary near the interface $\Sigma$ in order to resolve the highly
oscillatory SPP. In this subsection we give a brief overview of a local
mesh adaptation strategy based on \emph{goal-oriented} a posteriori error
estimation. This leads to a substantial saving in computational costs (see
Section~\ref{sec:computationalresults}), due to local refinement, while
maintaining an optimal convergence order in a quantity of
interest\cite{becker2001}. We focus primarily on aspects of
implementation---a full, detailed discussion is beyond the scope of this
paper and will be the subject of future work.

An efficient method for \emph{a posteriori} error control is the \emph{dual
weighted residual} (DWR) method developed by Becker and
Rannacher~\cite{becker1996a,becker1996b,becker2001}. This method constructs
estimates of local error contributions in terms of a target functional $\J$
with the help of a ``dual problem''. More precisely, let $\J(\vE)$ be a
quantity of interest given by a possibly non-linear Gâteaux-differentiable
function, viz.,
\begin{align}
  \J\,:\,\vec H(\text{curl};\Omega) \to\C.
\end{align}
The corresponding dual problem to (\ref{eq:variationalformulation}) is to
find a solution $\vZ\in\vec H(\text{curl};\Omega)$ such that
\begin{align}
  A(\vp,\vZ) =\;\; \text{D}_{\vE}\mathcal{J}(\vE)[\vp]
  \label{eq:dualproblem}
\end{align}
for all $\vp\in\vec H(\text{curl};\Omega)$. Here, $\text{D}_{\vE}.[\vp]$
denotes the Gâteaux derivative in direction $\vp$ with respect to $\vE$.
The following result can be directly applied to variational
problem~\eqref{eq:variationalformulation}.
\medskip

\begin{proposition}[Becker and Rannacher~\cite{becker2001}]
  Let $\vE$ and $\vZ$ be the solution of (\ref{eq:variationalformulation})
  and (\ref{eq:dualproblem}), respectively. Let $\vE_H$ and $\vZ_H$ be
  finite element approximations of the primal and dual solution associated
  with a discretization $\TTH(\Omega)$ of $\Omega$. Then, up to a term $R$
  of higher order:
  \begin{gather}
    \big|\mathcal{J}(\vE)-\mathcal{J}(\vU)\big|
    \le \sum_{Q\in\,\TTH}\eta_Q
    +R,\quad\text{with}
    \\
    \eta_Q:=
    \frac 12\,\Big|
    \rho_Q(\vE_H,\vZ-\vZ_H) +
    \rho^\ast_Q(\vZ_H,\vE-\vE_H)
    \Big|.
    \label{eq:localindicator}
  \end{gather}
\end{proposition}
\medskip

Here, $\rho_Q$ and $\rho^\ast_Q$ denote the primal and dual cell-wise
residual, respectively, associated with variational equations
(\ref{eq:variationalformulation}) and (\ref{eq:dualproblem}):
\begin{align}
  \label{eq:rhoQ}
  \rho_Q(\vE_H,\vZ-\vZ_H) &=
  F\big((\vZ-\vZ_H)\chi_Q\big) - A\big(\vE_H,(\vZ-\vZ_H)\chi_Q\big),
  \\[0.1em]
  \label{eq:rhoastQ}
  \rho^\ast_Q(\vZ_H,\vE-\vE_H) &=
  \text{D}_{\vE}\mathcal{J}(\vE)[(\vE-\vE_H)\chi_Q] -
  A\big((\vE-\vE_H)\chi_Q,\vE_H\big),
\end{align}
with the indicator function $\chi_Q$ that is 1 on the cell $Q$, and 0
otherwise. The \emph{local error indicator} $\eta_Q$ given by
\eqref{eq:localindicator} can now be approximated and used in a
local refinement strategy \cite{becker2001}.

Our goal is an optimal local refinement for the numerical observation of
SPPs on the conducting sheet $\Sigma$. In principle, a number of choices
for the quantity of interest, $\mathcal{J}(\vE)$, are possible in order to
achieve this goal. Here, we choose the quantity
\begin{align}
  \mathcal{J}(\vE):=\int_\Omega\varpi(\vx)\,\big\|\nabla\times\vE\,\big\|^2
  \rm{d}\vx,
  \label{eq:quantityofinterest}
\end{align}
with an appropriate (essentially bounded), non-negative weighting function
$\varpi(\vx)$ that localizes the integral around the interface $\Sigma$;
see Section~\ref{sec:computationalresults} for the concrete choice of
$\varpi$ for our simulations.
Choice \eqref{eq:quantityofinterest} for the quantity of interest leads to
a localized right-hand side $\mathcal{J}$ of the dual problem that is
sensitive to the highly oscillatory SPP  associated with the electric
field, $E$. Consequently, the \emph{weight} $\vZ-\vZ_H$ in residual
\eqref{eq:rhoQ} is generally large near the interface and at points where
the influence of the solution on quantity \eqref{eq:quantityofinterest} is
high.
\begin{remark}
  A purely residual-based error estimator on the other hand corresponds to
  a uniform weight distribution. The DWR method (with our choice of right
  hand side) will lead to a more localized refinement.
\end{remark}

In order to guarantee a uniform refinement over the interface $\Sigma$, the
local refinement strategy is augmented by \emph{additionally} selecting all
cells $Q$ for refinement that fulfill
\begin{align}
  1-(1/2)^{\text{\#cycle}-1}
  \;\le\;
  \varpi(x_Q)/\max(\varpi),
\end{align}
where $x_Q$ denotes the center of $Q$.

\medskip

\begin{remark}[Evaluation of residuals]
  A classical approach to evaluate (\ref{eq:rhoQ}) and (\ref{eq:rhoastQ})
  is to use a higher-order approximation for $\vZ$ and $\vE$ and transform
  $\rho_Q$ and $\rho^\ast_Q$ into a strong residual form~\cite{becker2001}.
  We follow a different strategy that does not involve solving higher-order
  solutions. Instead of a higher-order approximation of $\vZ,$ we use a
  patch-wise projection $\pi_{2H}^{(2)}\vZ_H$ to a higher-order space on a
  coarser mesh level~\cite{richter2006}, viz.,
  \begin{align}\label{eq:projection-est}
    \vZ-\vZ_H \approx \pi_{2H}^{(2)}\vZ_H - \vZ_H,
    \quad
    \vE-\vE_H \approx \pi_{2H}^{(2)}\vE_H - \vE_H,
  \end{align}
  in combination with the (variational) residuals \eqref{eq:rhoQ} and
  \eqref{eq:rhoastQ}.
\end{remark}

\medskip

\begin{remark}[Mixed form of error estimator]
  Error estimator~\eqref{eq:projection-est} also has a form with only the
  primal residual $\rho(\,.\,)(\,.\,)$ appearing. However, the mixed form
  of the above error estimator should be used here to ensure adequate,
  simultaneous refinement not just only with respect to a localized SPP (on
  the interface), but also with respect to a singular right-hand side
  $\vJa$ (modeling of a Hertzian dipole).
\end{remark}


\section{Numerical results}
\label{sec:computationalresults}

In this section, we present a number of computational results pertaining to
the excitation by a vertical electric dipole of an SPP on a 2D conducting
sheet. First, we construct a PML which involves use of a certain parameter,
$s_0$. Second, we demonstrate the necessity for a PML and carry out a study
for the suitable choice of $s_0$. Third, we compare the analytical results
of Section~\ref{sec:analytical} to (direct) numerical simulations based on
our finite element formulation of Section~\ref{sec:numerics}. All numerical
computations are carried out with the \texttt{C++} finite element toolkit
\texttt{deal.II}~\cite{dealii82,dealii83}. We use the direct solver
\texttt{UMFPACK}~\cite{davis11,suitesparse4} for the resulting linear
system of equations.

\subsection{Setup}

In order to carry out the numerical simulations, we consider a vertical
electric dipole positioned at $\vec a = (0,0.75)$ and $\vec a = (0,1)$
above the conducting sheet. Recall that the corresponding current densities
are
\begin{align}
  \vJa=\begin{pmatrix}0\\J_0\end{pmatrix}\,\delta(\vx - \vec a),
  \qquad
  \vMa=0.
\end{align}
The conducting sheet has surface conductivity $\ssr=(\sigma_1+i\sigma_2)I$
in tensor form ($I$: unit second-rank tensor). Below, we carry out a
parameter study with different values of constant scalar $\ssr$.

The computational domain $\Omega$ is chosen to be the ball with radius
$R=8\pi$, which is 8 times the free-space wavelength of the dipole
radiation. A PML is used in the outer region at distance
$\rho>0.8\,R$ from the origin; see Section~\ref{sec:pml}. Following
Monk~\cite{monk03}, the (nonnegative) scaling function $s(\rho)$ is chosen
to be
\begin{align}
  s(\rho) = s_0 \frac{(\rho - 0.8 R)^2}{(R-0.8 R)^2},
\end{align}
where $s_0$ is a free parameter; the suitable choice of $s_0$ is discussed
in Section~\ref{sub:pollutionandpml}.

We use the adaptive refinement strategy that was outlined in
Section~\ref{sec:aposteriori} with the quantity of interest
(\ref{eq:quantityofinterest}) and the choice
\begin{align}
  \varpi(\vx) =
  \begin{cases}
    \begin{aligned}
      &\cos^2\Big(\frac\pi2\,\frac{y}{d_\varpi}\Big),
      &\quad&\text{if }|y|\le d_\varpi,
      \\[0.2em]
      &0 &\quad&\text{otherwise}.
      \quad
    \end{aligned}
  \end{cases}
  \label{eq:interfacereg}
\end{align}
The parameter $d_\varpi$ controls the width of the region around the
interface that is uniformly refined and is chosen to be $d_\varpi=1.5625$,
a few multiples of the SPP wavelength $2\pi/(\text{Re}\,\kmr)$. The current
density of the dipole source is regularized according to
\begin{align}
  \delta(\vx-\vec a)\approx
  \begin{cases}
    \begin{aligned}
      &\Big(\frac{\pi}{2}-\frac{2}{\pi}\Big)^{-1} d^{-2}
      \cos^2\Big(\frac\pi2\,\frac{\|\vx - \vec a\|}{d}\Big)
      & &\text{for }\|\vx-\vec a\|<d,
      \\[0.2em]
      &0 & &\text{otherwise,}
    \end{aligned}
  \end{cases}
  \label{eq:dipolereg}
\end{align}
with $d=0.15625$. The regularization parameter $d$ should be chosen as
small as possible and has to be scaled with the initial mesh size such that
$\vJa$ is always well integrated numerically.
Figure~\ref{fig:graphenesheet} summarizes the aforementioned setup.


\subsection{On perfectly matched layer}
\label{sub:pollutionandpml}

\begin{figure}[!t]
  \centering
  \subfloat[Scattered $\text{Re}\,(\vE_x)$ as function of position $x$]{
    \includegraphics{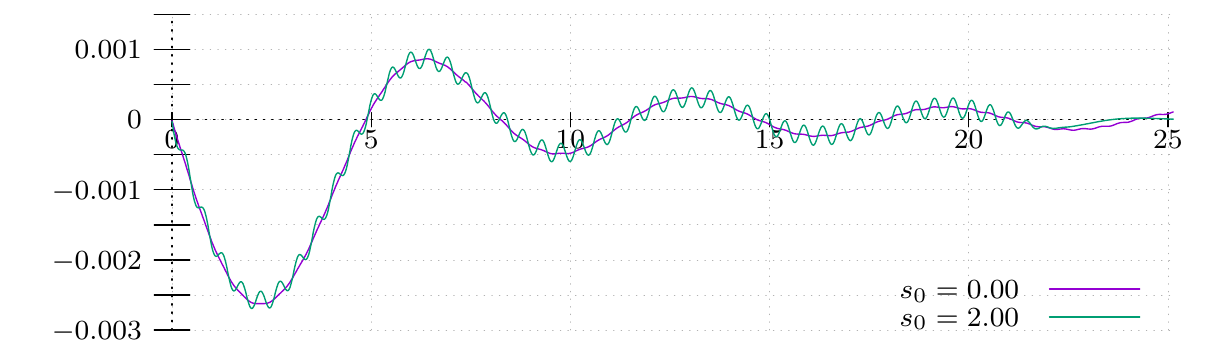}
  }

  \subfloat[zoom into the PML, $s_0=0.0$, $0.25$, $0.5$, $1.0$]{
    \includegraphics{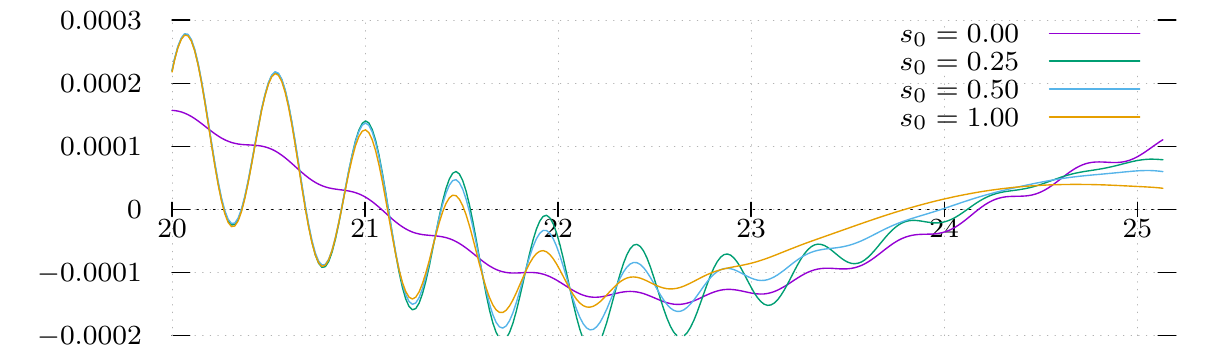}
  }

  \subfloat[zoom into the PML, $s_0=1.0$, $2.0$, $4.0$, $8.0$]{
    \includegraphics{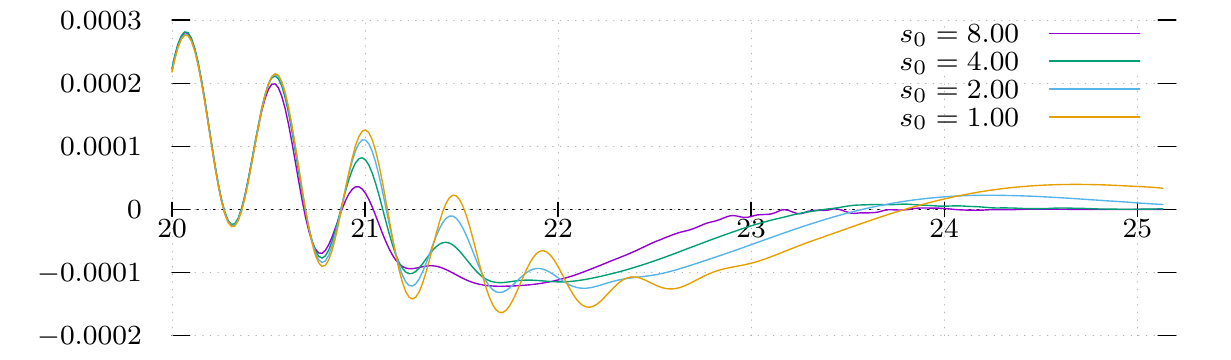}
  }

  \caption{
    Real part of $x$-component of scattered field, $\text{Re}E_x^{{\rm
    sc}}=\text{Re}(E_x-E_x^{{\rm pr}})$, on the interface $\Sigma$ as a
    function of position, $x$; $E_x^{pr}$ is the primary dipole field
    component, in the absence of the sheet.  (a) Comparison of no PML with
    a PML of parameter $s_0=2$; (b) zoom into the PML with parameter values
    $s_0=0, 0.25, 0.5, 1.0$; (c) zoom into the PML with $s_0=1.0, 2.0, 4.0,
    8.0$.}
  \label{fig:pmltest}
\end{figure}

In this subsection, we demonstrate the necessity for a PML. In our
numerical setup, the challenging part of a direct finite element simulation
is the two-scale character that the electromagnetic wave exhibits in the
spatial resolution. Recall that the desired SPP has a wavelength much
smaller than the one manifested by the dipole free-space radiation field.

In particular, we are interested in observing SPPs with an associated wave
number $\text{Re}\,\kmr\approx 10-100$, compared to the (rescaled)
free-space wave number $k_r=1$. Moreover, the amplitude of the SPP scales
exponentially with the distance, $a$, of the dipole from the interface;
hence, certain configurations of physical appeal exhibit a ratio in
amplitude of about 1:10 to 1:1000 between the SPPs and the dipole
free-space radiation field. It turns out that absorbing boundary condition
\eqref{eq:absorbingbc} that we use--although it is a first-order absorbing
boundary condition--is \emph{not} well suited for the numerical study of
the SPP. A significant suppression of the SPP amplitude can be evident.

In order to study the influence of the absorbing boundary condition on the
SPP, we perform a series of numerical simulations of the geometry given in
Figure~\ref{fig:graphenesheet}, for different values of $s_0$ which
controls the absorption strength of the PML; specifically, we use $s_0=0$,
$0.25$, $0.5$, $1.0$, $2.0$, $4.0$, and $8.0$, respectively. The material
parameters are set to $\mur=\er=I$. We choose $\vec a=(0,1)$ for the
position vector of the dipole, with sheet of relatively small, purely
imaginary  surface conductivity, i.e.,
\begin{align}
  \ssr = i\sigma_2 = 0.15\,i.
\end{align}
This value corresponds to a relatively weak SPP, i.e., with a relatively
small amplitude, that has a purely real-valued wave number, $\kmr\approx
13.33$. This can be understood as follows: Because
$\text{Re}\sigma^\Sigma_r=0$, the sheet does not cause any dissipation and,
thus, the corresponding SPP does not exhibit any decay.

In order to examine the influence of different boundary conditions and
choices of parameters in some detail, we extract the $x$-component of the
scattered electric field, $E_x^{\rm sc}=E_x-E_x^{\rm pr}$, on the interface
$\Sigma$, and plot this component as a function of $x$, where $\vE^{\rm
pr}$ is the (primary) dipole field in the absence of the sheet; see
Figure~\ref{fig:pmltest}. Boundary condition \eqref{eq:absorbingbc} without
PML ($s_0=0$) has a strong influence on the observed scattered field. While
branch-cut contribution (\ref{eq:branchcutcontributionrescaled}) of the
scattered field can be observed in the numerical simulation, pole
contribution (\ref{eq:polecontributionrescaled}), which is responsible for
the fast oscillation with wave number $\kmr$, is suppressed
(Figure~\ref{fig:pmltest}a). This property can be explained by the fact
that boundary condition (\ref{eq:absorbingbc}), viz.,
\begin{align*}
  \vec n\times\vB+\sqrt{\muri\er}\vE_T = 0,
\end{align*}
is only applicable to dipole radiation with a wave number $k_r=1$. In the
case where the pole contribution is characterized, e.g., by
$\kmr\approx13.3$, this boundary condition causes a reflection that results
in a suppression of the fast spatial oscillation of the SPP. In light of
the parameter study (Fig.~\ref{fig:pmltest}), we choose a PML with strength
$s_0=2$ for all subsequent computations. This is a balanced choice between
a PML that is strong enough to minimize the influence of the boundary
($s_0\ge2$) and the unwanted influence of the PML ($s_0\le2$);
cf.~Figure~\ref{fig:pmltest}c.


\subsection{Comparison of numerical results to analytical solution}
\label{sub:comparison}

\begin{figure}[!t]
  \centering
  \subfloat[weak absorption]{
    \mbox{\hspace{1.43cm}}
    \includegraphics[height=6cm]{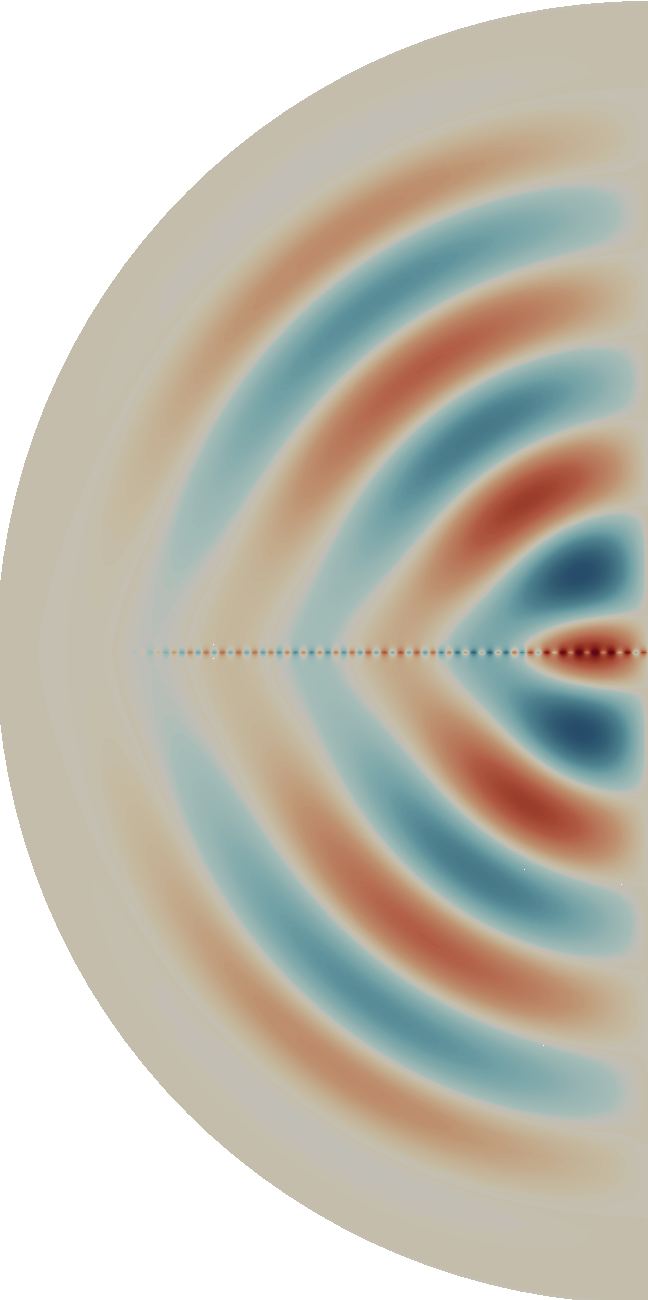}
  }
  \subfloat[strong absorption]{
    \includegraphics[height=6cm]{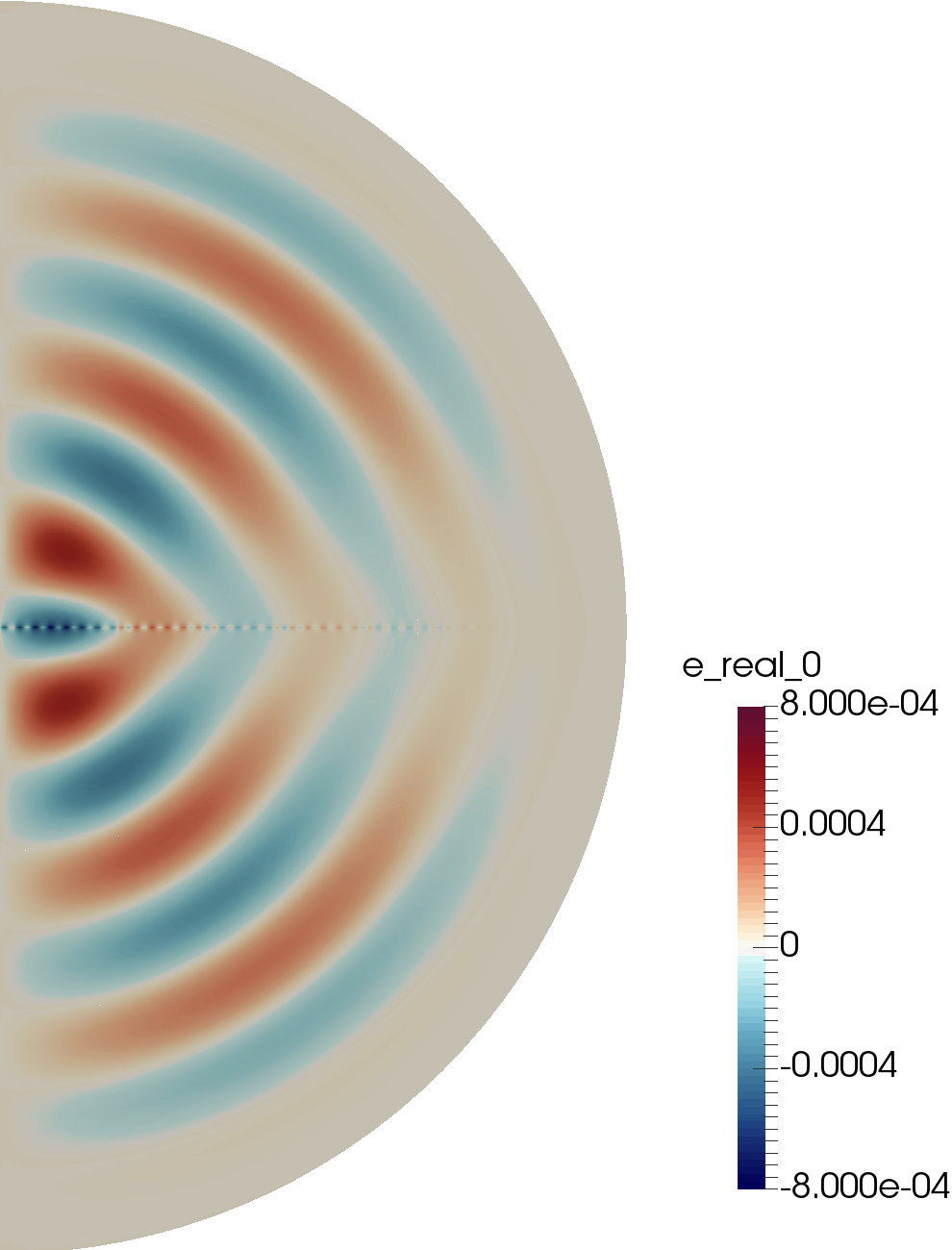}
  }

  \subfloat[Zoom of (b) around origin with detailed SPP]{
    \includegraphics[width=10cm]{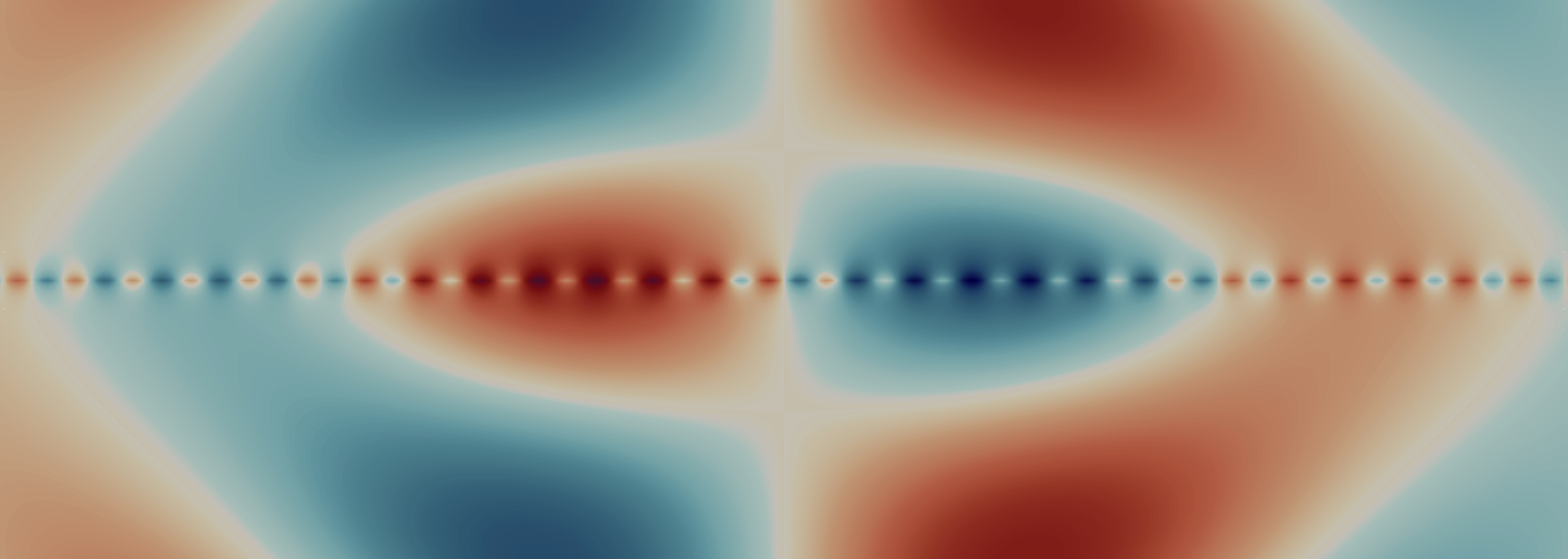}
  }
  \caption{
    Real part of $x$-directed scattered electric field,
    $\text{Re}\,(E_x^{\rm sc})$, on interface $\Sigma$, computed for: (a)
    $\ssr=4.0\times10^{-4}+0.2i$ ($\km\approx10.0+0.02i$); and (b)
    $\ssr=2.0\times10^{-3}+0.2i$ ($\km\approx10.0+0.1i$). The dipole
    elevation distance is $a=1.00$.}
  \label{fig:scatteredfield}
\end{figure}

\begin{figure}[!htbp]
  \centering
  \subfloat[full mesh]{
    \includegraphics[height=4.0cm]{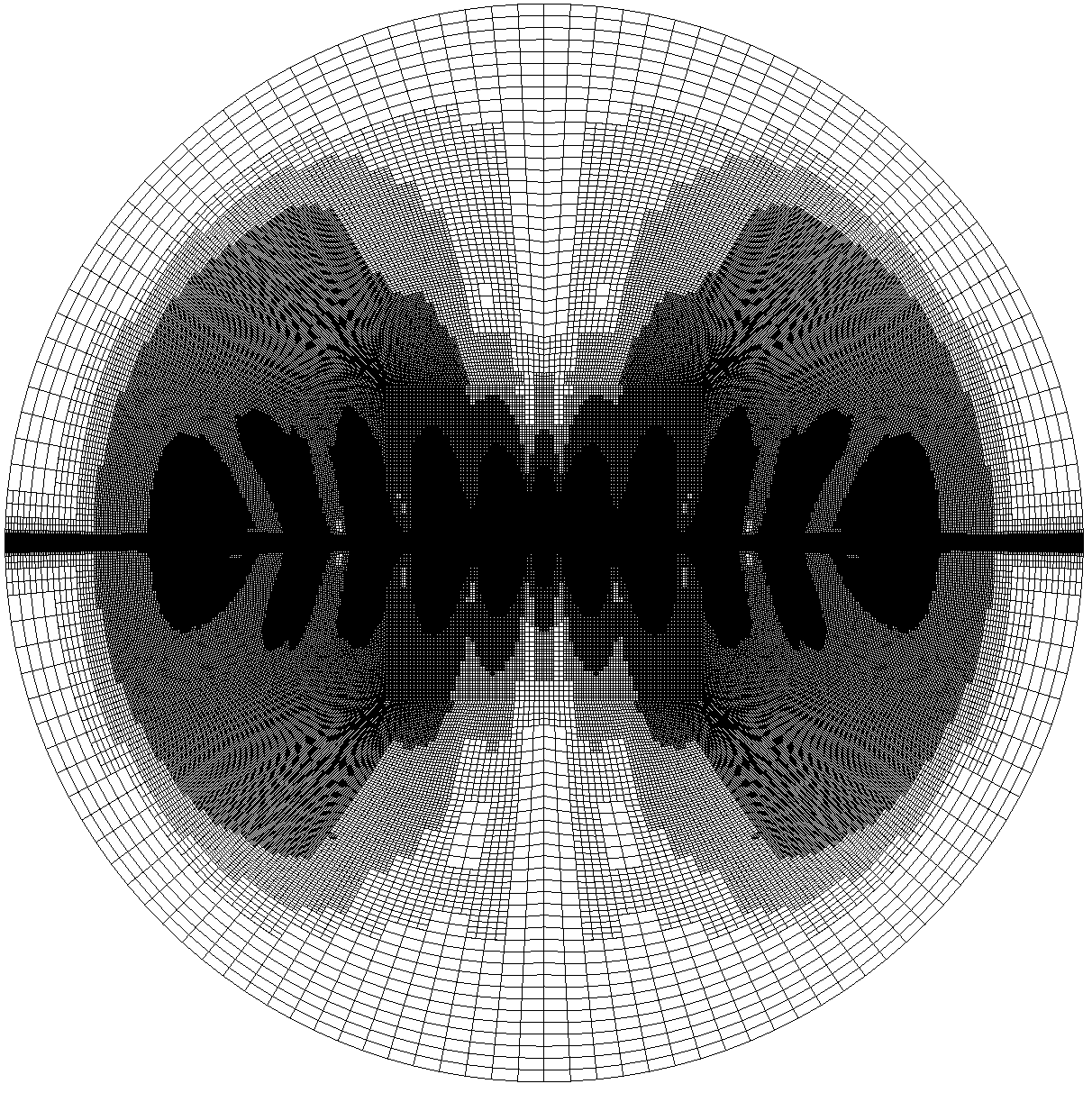}
  }
  \subfloat[zoom]{
    \includegraphics[height=4.0cm]{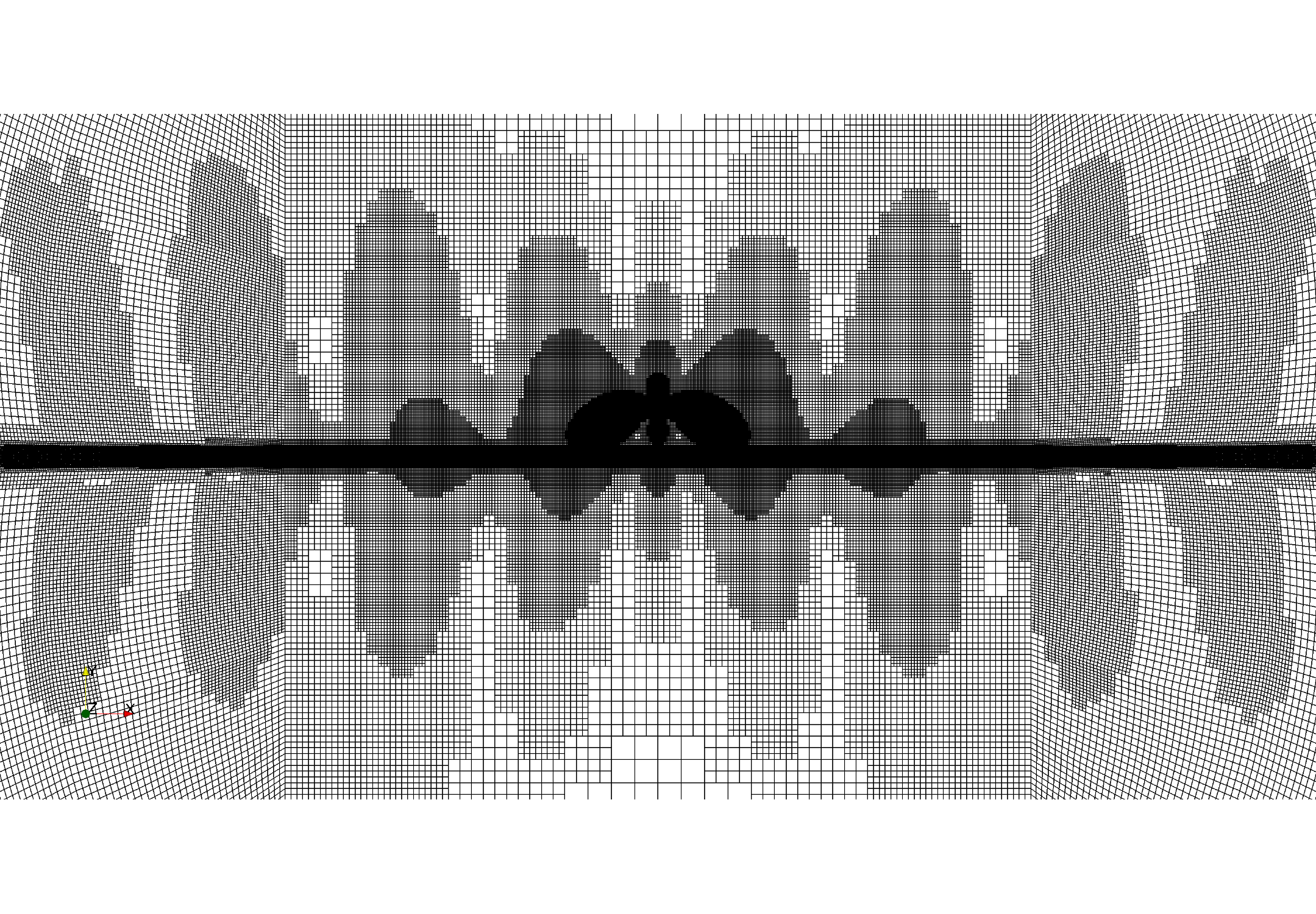}
  }
  \caption{
    The locally refined mesh (a) obtained with the adaptive method outlined in
    Section~\ref{sec:aposteriori}. The mesh has a total number of around
    $200$~thousand cells; the finest resolution around the interface
    corresponds to a uniformly refined mesh of around $5$~million cells.}
  \label{fig:localrefinement}
\end{figure}

In this subsection, we compare our numerical results to the analytical
solution of Section~\ref{sec:analytical}. In particular, we expect to
observe the SPP described by~(\ref{eq:polecontributionrescaled}). The
(complex-valued) wave number $\kmr$ associated with this SPP scales with
the surface conductivity $\ssr$ as follows (given in rescaled quantities,
as explained in~Section~\ref{sec:rescaling}):
\begin{align}
  \kmr =
  \sqrt{\mur\er-\frac{4\mur^2\er^2}{(\ssr)^2}};
  \qquad
  \kmr \approx
  \frac{2i\,\mur\er}{\sigma_r^\Sigma}
  \quad
  \text{if}\quad |2\sqrt{\mur\er}|\gg|\sigma_r^\Sigma|.
  \label{eq:polaritonscalingrescaled}
\end{align}
In order to test our numerical method (Section~\ref{sec:numerics}) against
the analytical results (Section~\ref{sec:analytical}), we carry out a
parameter study for a variety of different values of $\ssr$; see
Table~\ref{tab:sigma}.
\begin{table}[!t]
  \center
  \footnotesize
  \begin{tabular} {c c}
    \toprule
    Surface conductivity $\ssr$ & Predicted $\kmr$
    \\[0.5em]
    %
    $2.56\times10^{-4}+0.160i$ &  $12.5+0.02i$  \\
    $1.78\times10^{-4}+0.133i$ &  $15.0+0.02i$  \\[0.3em]
    $1.28\times10^{-3}+0.160i$ &   $12.5+0.1i$  \\
    $8.89\times10^{-4}+0.133i$ &   $15.0+0.1i$  \\
    \bottomrule
  \end{tabular}
  \caption{Values of $\ssr$ used in the parameter study along with the predicted SPP
    wave numbers $\kmr$ by (\ref{eq:polaritonscalingrescaled}).
  }
  \label{tab:sigma}
\end{table}

Our computations are performed for dipole elevation distances $a=0.75,
1.00$. For each choice of parameters, we start with a relatively coarse
mesh (with 10k degrees of freedom) for the numerical simulation and run 6
mesh-adaptation cycles (with approximately 1.6M degrees of freedom on the
finest mesh.) Figure~\ref{fig:scatteredfield} shows the $x$-directed
scattered field for $\ssr=4.0\times10^{-4}+0.2i$ ($\km\approx10.0+0.02i$),
and $\ssr=2.0\times10^{-3}+0.2i$ ($\km\approx10.0+0.1i$). For this choice
of parameters, a relatively strong pole contribution can be observed as
shown in Figure~\ref{fig:scatteredfield}(b).
Figure~\ref{fig:localrefinement} shows the locally refined mesh after the
final refinement cycle.

For a comparison of our numerical method to the analytical results of
Section~\ref{sec:analytical}, the pole contribution
(\ref{eq:polecontributionrescaled}) and branch-cut contribution
(\ref{eq:branchcutcontributionrescaled}) are computed numerically. For this
purpose, we use a summed trapezoidal rule to evaluate the two integrals
involved in the branch-cut contribution. For the improper integral over
$(0, \infty)$, we further exploit an exponential decay of the integrand (as
a function of the integration variable, $\varsigma$) and introduce a cutoff
at $\varsigma\approx 1/\sqrt{hx}$, where $h$ is the interval size of the
summed trapezoidal rule. This $h$ is chosen adaptively in an iterative
cycle such that the relative change in the value of the integral between
$2h$ and $h$ is less than $0.5\,\%$.

\begin{figure}[!pt]
  \centering

  \subfloat[]{
    \includegraphics{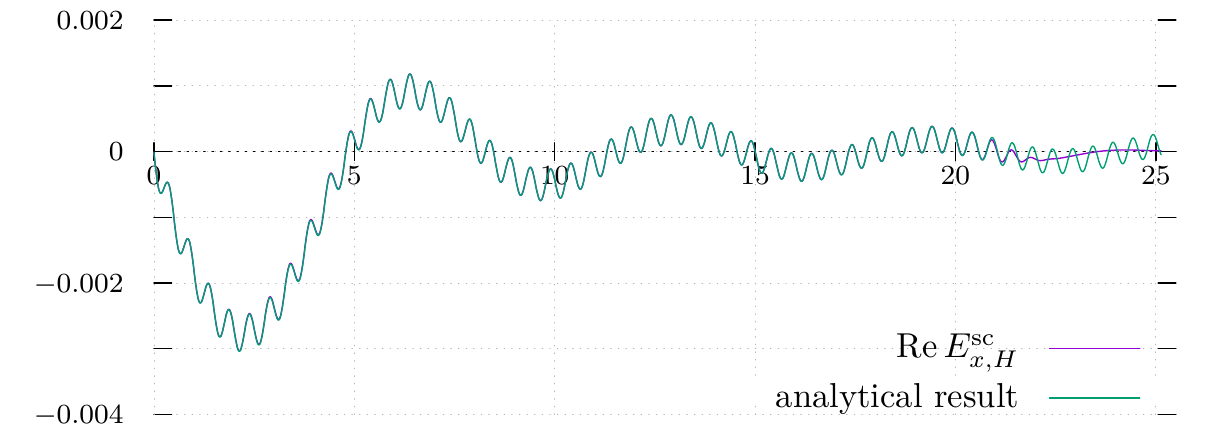}
  }

  \subfloat[]{
    \includegraphics{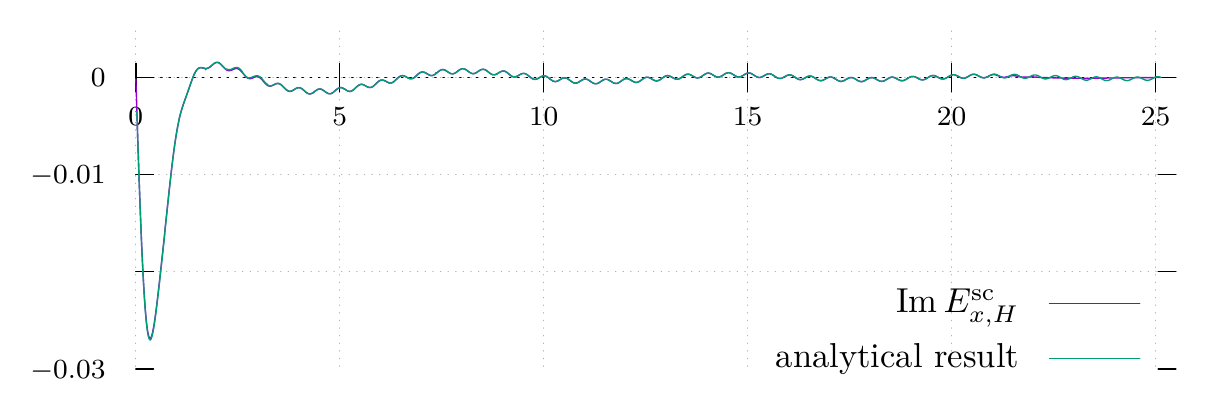}
  }

  \caption{
    Real part (a) and imaginary part (b) of the $x$-component of scattered
    electric field, $E_x^{\rm sc}$, as a function of position $x$ on
    interface $\Sigma$, in the presence of a vertical electric dipole at
    distance $a$ from the conducting sheet. The plots (a-b) show the
    numerical simulations based on our method as well as the analytical
    results computed by (\ref{eq:polecontributionrescaled}) and
    (\ref{eq:branchcutcontributionrescaled}) for the values $a=1.00$,
    $\kmr=12.5+0.02 i$.}
  \label{fig:parameterstudyinterface}
\end{figure}

\begin{figure}[!pt]
  \centering
  \subfloat[]{
    \includegraphics{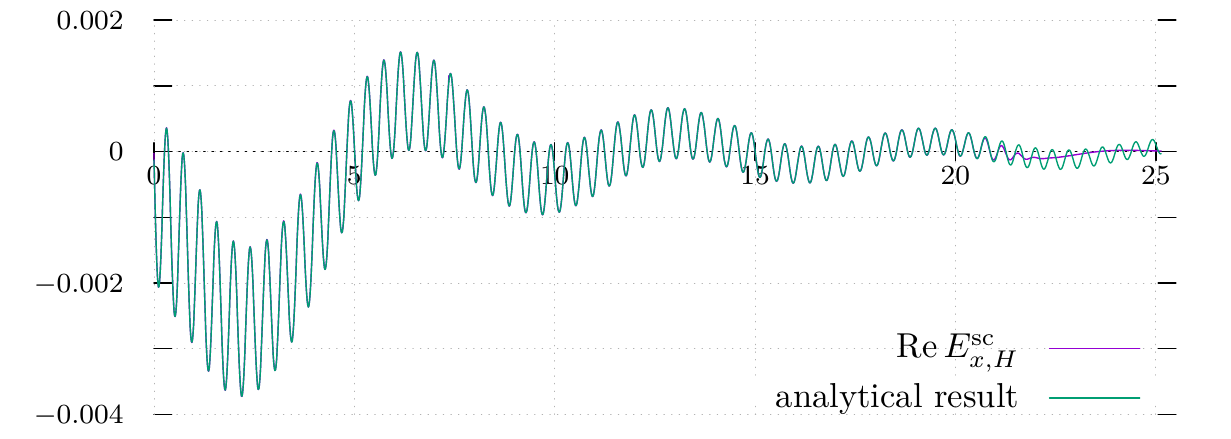}
  }

  \subfloat[]{
    \includegraphics{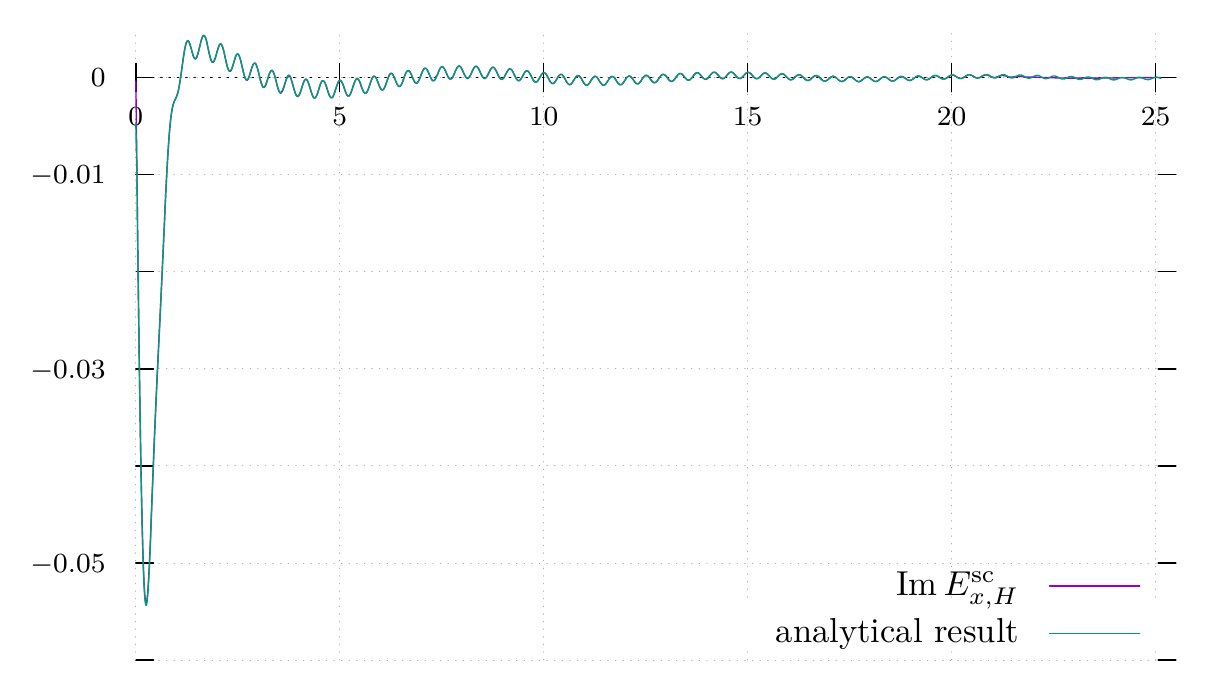}
  }
  \caption{
    Real part (a) and imaginary part (b) of the $x$-component of scattered
    electric field, $E_x^{\rm sc}$, as a function of position $x$ on
    interface $\Sigma$, in the presence of a vertical electric dipole at
    distance $a$ from the conducting sheet. The plots (a-b) show the
    numerical simulations based on our method as well as the analytical
    results computed by (\ref{eq:polecontributionrescaled}) and
    (\ref{eq:branchcutcontributionrescaled}) for the values $a=0.75$,
    $\kmr=15.0+0.1 i$.}
  \label{fig:parameterstudyinterface2}
\end{figure}

In Figures~\ref{fig:parameterstudyinterface} and
\ref{fig:parameterstudyinterface2}, the analytical and numerical results
are compared graphically for dipole elevation distances $a=0.75, 1.00$ and
$\kmr=12.5+0.02 i, 15.0+0.1 i$. The real and imaginary parts of the
scattered electric field in the $x$-direction, $E_x^{\rm sc}$, are plotted
as a function of position $x$ on the interface, $\Sigma$. It is evident
that the numerical and analytical results are in excellent agreement
outside the PML ($0\le x\le 20$).

Table~\ref{tab:comparison} summarizes the parameter study quantitatively.
The $L^2$-error between numerical and analytical result is computed outside
of the PML on the interface $\Sigma$. The SPP is generally well
approximated after the 4th cycle with a resolution of around 100 thousand
degrees of freedom. According to \eqref{eq:convergenceorder}, we expect a
convergence order of
$\|\vE_{h,T}-\vE_T\|_{L^2(\Omega)}\sim\mathcal{O}(\#\text{dofs})$. Indeed,
we observe a linear convergence of the error with respect to the number of
refinement steps, and thus, a linear convergence in number of degrees of
freedom.

\begin{table}[p]
  \center
  \subfloat[$a=0.75$, $k_m\approx12.5+0.02i$, $k_m\approx12.5+0.1i$]{
    \footnotesize
    \begin{tabular} {r r r c c r r c c}
      \toprule
      & \multicolumn{4}{c}{$\kmr\approx12.5+0.02i$} &
      \multicolumn{4}{c}{$\kmr\approx12.5+0.1i$}\\
      \cmidrule(lr){2-5}
      \cmidrule(lr){6-9}
      Cycle  & Cells  & DoFs    & $L^2$-error & & Cells  & DoFs    & $L^2$-error &
      \\[0.5em]
      1 &   1280 &   10304 & 2.49e-2 &     - &   1280 &   10304 & 1.52e-2 &     - \\
      2 &   2624 &   21660 & 7.08e-1 & -4.83 &   2624 &   21660 & 2.03e+0 & -7.06 \\
      3 &   5474 &   45580 & 2.48e-2 &  4.84 &   5450 &   45388 & 1.41e-2 &  7.17 \\
      4 &  11930 &   99020 & 1.41e-2 &  0.82 &  11882 &   98636 & 8.09e-3 &  0.80 \\
      5 &  27908 &  229872 & 7.00e-3 &  1.01 &  27842 &  229344 & 4.03e-3 &  1.00 \\
      6 &  71246 &  582828 & 3.60e-3 &  0.96 &  71114 &  581796 & 2.03e-3 &  0.99 \\
      7 & 191906 & 1561336 & 1.81e-3 &  0.99 & 191672 & 1559404 & 1.02e-3 &  1.00 \\
      \bottomrule
    \end{tabular}
  }

  \subfloat[$a=0.75$, $k_m\approx15.0+0.02i$, $k_m\approx15.0+0.1i$]{
    \footnotesize
    \begin{tabular} {r r r c c r r c c}
      \toprule
      & \multicolumn{4}{c}{$\kmr\approx15.0+0.02i$} &
      \multicolumn{4}{c}{$\kmr\approx15.0+0.1i$}\\
      \cmidrule(lr){2-5}
      \cmidrule(lr){6-9}
      Cycle  & Cells  & DoFs    & $L^2$-error & & Cells  & DoFs    & $L^2$-error &
      \\[0.5em]
      1 &   1280 &   10304 & 5.33e-3 &     - &   1280 &   10304 & 4.35e-3 &     - \\
      2 &   2624 &   21660 & 4.94e-3 &  0.11 &   2624 &   21660 & 1.24e-2 & -1.51 \\
      3 &   5474 &   45556 & 9.10e-3 & -0.88 &   5474 &   45556 & 4.40e-3 &  1.50 \\
      4 &  11906 &   98828 & 4.17e-3 &  1.13 &  11906 &   98828 & 2.18e-3 &  1.01 \\
      5 &  27890 &  229704 & 1.78e-3 &  1.23 &  27890 &  229704 & 1.04e-3 &  1.06 \\
      6 &  71228 &  582684 & 9.28e-4 &  0.94 &  71228 &  582672 & 5.27e-4 &  0.99 \\
      7 & 191978 & 1561852 & 4.69e-4 &  0.99 & 191978 & 1561852 & 2.64e-4 &  1.00 \\
      \bottomrule
    \end{tabular}
  }

  \subfloat[$a=1.00$, $k_m\approx12.5+0.02i$, $k_m\approx12.5+0.1i$]{
    \footnotesize
    \begin{tabular} {r r r c c r r c c}
      \toprule
      & \multicolumn{4}{c}{$\kmr\approx12.5+0.02i$} &
      \multicolumn{4}{c}{$\kmr\approx12.5+0.1i$}\\
      \cmidrule(lr){2-5}
      \cmidrule(lr){6-9}
      Cycle  & Cells  & DoFs    & $L^2$-error & & Cells  & DoFs    & $L^2$-error &
      \\[0.5em]
      1 &   1280 &   10304 & 2.96e-3 &     -  &   1280 &   10304 & 2.81e-3 &     - \\
      2 &   2600 &   21504 & 1.37e-1 & -5.54  &   2600 &   21504 & 3.94e-1 & -7.13 \\
      3 &   5402 &   44980 & 1.33e-3 &  6.69  &   5402 &   44980 & 8.96e-4 &  8.78 \\
      4 &  11786 &   97868 & 6.39e-4 &  1.06  &  11786 &   97868 & 4.02e-4 &  1.16 \\
      5 &  27728 &  228240 & 3.26e-4 &  0.97  &  27728 &  228240 & 2.06e-4 &  0.96 \\
      6 &  71366 &  584010 & 1.68e-4 &  0.96  &  71366 &  583998 & 1.04e-4 &  0.99 \\
      7 & 194786 & 1584772 & 8.74e-5 &  0.94  & 194786 & 1584760 & 5.91e-5 &  0.82 \\
      \bottomrule
    \end{tabular}
  }

  \subfloat[$a=1.00$, $k_m\approx15.0+0.02i$, $k_m\approx15.0+0.1i$]{
    \footnotesize
    \begin{tabular} {r r r c c r r c c}
      \toprule
      & \multicolumn{4}{c}{$\kmr\approx15.0+0.02i$} &
      \multicolumn{4}{c}{$\kmr\approx15.0+0.1i$}\\
      \cmidrule(lr){2-5}
      \cmidrule(lr){6-9}
      Cycle  & Cells  & DoFs    & $L^2$-error & & Cells  & DoFs    & $L^2$-error &
      \\[0.5em]
      1 &   1280 &   10304 & 2.30e-3 &    -  &   1280 &   10304 & 2.20e-3 &    - \\
      2 &   2600 &   21504 & 6.43e-4 & 1.84  &   2600 &   21504 & 2.20e-3 & 0.00 \\
      3 &   5426 &   45196 & 6.24e-4 & 0.04  &   5426 &   45196 & 4.38e-4 & 2.32 \\
      4 &  11834 &   98276 & 1.84e-4 & 1.76  &  11834 &   98276 & 1.66e-4 & 1.40 \\
      5 &  27812 &  228876 & 9.29e-5 & 0.98  &  27812 &  228876 & 8.67e-5 & 0.94 \\
      6 &  71546 &  585474 & 4.96e-5 & 0.90  &  71546 &  585474 & 4.60e-5 & 0.91 \\
      7 & 195095 & 1587232 & 3.30e-5 & 0.59  & 195095 & 1587232 & 3.15e-5 & 0.55 \\
      \bottomrule
    \end{tabular}
  }

  \caption{Convergence history and $L^2$ error of the real part of the
    scattered electric field in $x$-direction between the numerical and
    analytical solution
      computed outside of the PML on the interface $\Sigma$.
    The columns after the $L^2$ error are the $\log_2$ reduction rates of
    the error. }
  \label{tab:comparison}
\end{table}


\section{Conclusion and outlook}
\label{sec:conclusion}

In this paper, we developed a variational framework for the numerical
simulation of electromagnetic SPPs excited by current-carrying sources
along an infinitely thin conducting sheet, e.g., single-layer graphene. The
sheet is modeled by an idealized, oriented hypersurface. The effect of the
induced surface current of the sheet was taken into account in the
corresponding boundary value problem of Maxwell's equations via jump
condition~\eqref{eq:jumpcondition} for the tangential component of the
magnetic field; this jump is proportional to the surface current density.
The conductivity of the sheet is a parameter that controls the
strength of this discontinuity.

One of the main advantages offered by our approach is the natural
incorporation of the jump condition in a variational formulation as a weak
discontinuity, \emph{without} regularization of the interface by a layer of
finite thickness. We tested our numerical treatment against analytical
predictions in the case with a vertical dipole radiating over an infinite
conducting sheet in 2D, and observed excellent agreement. In our numerical
treatment, a linear asymptotic reduction rate could be observed for all
testcases. Notably, the use of an adaptive local refinement procedure
within our approach  achieved significant economy in the total number of
degrees of freedom in comparison to uniform mesh refinement.

Our numerics admit several generalizations and extensions. For instance,
the treatment of the jump condition as a weak discontinuity over an
interface is \emph{not} limited to a simple (lower-dimensional) hyperplane;
it can be generalized to reasonably arbitrary hypersurfaces. In fact,
technically speaking, our variational framework can be readily used without
modifications to model any geometry that is meshable by quadrilaterals.
This flexibility should enable efficient numerical simulations of complex
geometries, e.g., waveguides that contain a few graphene
layers~\cite{yeh-book}. Although our numerical results
focused on 2D thus far, our underlying choice of local adaptivity can lead
to a significant reduction of computational cost in higher spatial
dimension.



\section*{References}

\bibliographystyle{elsarticle-num}

\bibliography{maiermargetisluskin}

\end{document}